\documentclass[reqno,10pt] {article}
\usepackage{graphicx, color,amsfonts,amsmath,dsfont,amssymb,hyperref,floatrow}
\usepackage{xcolor}
\usepackage{graphicx}
\usepackage{tikz}
 \usepackage{amsthm}

\newtheorem{theorem}{Theorem}
\newtheorem{remark}[theorem]{Remark}
\newtheorem{lemma}[theorem]{Lemma}
\newtheorem{proposition}[theorem]{Proposition}
\newtheorem{corollary}[theorem]{Corollary}
\newtheorem{definition}[theorem]{Definition}

\def \E{ \mathbb E  }

\definecolor{remi}{rgb}{1,0,0}
\usepackage{titlesec}
    \titleformat{\section}[hang]
        {\color{remi}{}\bfseries\filcenter\large}
        {\thesection.}
        {0.4em}
        {}[]
        
 \hypersetup{
    linktoc=page,
    linkcolor=blue,          
    citecolor=green,        
    filecolor=blue,      
    urlcolor=cyan,
    colorlinks=true           
}

\renewcommand{\tilde}{\widetilde}          
\DeclareMathSymbol{\leqslant}{\mathalpha}{AMSa}{"36} 
\DeclareMathSymbol{\geqslant}{\mathalpha}{AMSa}{"3E} 
\DeclareMathSymbol{\eset}{\mathalpha}{AMSb}{"3F}     
\renewcommand{\leq}{\;\leqslant\;}                   
\renewcommand{\geq}{\;\geqslant\;}                   
\newcommand{\p}{\mathbb{P}}
\newcommand{\R}{\mathbb{R}}

\newcommand{\N}{\mathbb{N}}

\newcommand{\ind}{\mathbf{1}}

\begin{document}

\title{Lognormal $\star$-scale invariant random measures}
\author{Romain Allez, R\'emi Rhodes \& Vincent Vargas\\
{\small Universit{\'e} Paris-Dauphine, Ceremade, UMR 7534,} \\{\small Place du marŽchal de Lattre de Tassigny, 75775 Paris Cedex 16, France.}}
\maketitle
\begin{abstract}
In this article, we consider the continuous analog of the celebrated Mandelbrot star equation with lognormal weights. Mandelbrot introduced this equation to characterize the law of multiplicative cascades. We show existence and uniqueness of measures satisfying the aforementioned continuous equation; these measures fall under the scope of the Gaussian multiplicative chaos theory developed by J.P. Kahane in 1985 (or possibly extensions of this theory). As a by product, we also obtain an explicit characterization of the covariance structure of these measures.  We also prove that qualitative properties such as long-range independence or isotropy can be read off  the  equation.
\end{abstract}
\begin{center}
\end{center}
\tableofcontents

\normalsize

\section{Introduction}
Fractality and the related concept of scale invariance is nowadays well introduced in many fields of applications ranging from physics, finance, information or social sciences. The scale-invariance property of
a stochastic process is usually quantified by the scaling exponents $\xi(q)$ associated with the power-law
behavior of the order $q$ moments of the �fluctuations� at different scales. More precisely, if $X_t$ is a $1$-d process with stationary increments, 
we can consider the $q$-th moments of its fluctuations at scale $l$:
$$\E\big[|X_{t+l}-X_t|^q\big].$$
The scaling exponents $\xi(q)$ are defined through the following  power-law scaling:
$$\E\big[|X_{t+l}-X_t|^q\big]=C_ql^{\xi(q)}\quad \forall l<T.$$

When $\xi(q)=qH$ is linear, the process is said to be monofractal. Famous examples of such processes are (fractional) Brownian motion, $\alpha$-stable L\'evy processes or Hermitte processes. When $\xi$ is nonlinear, the process is said to be multifractal. The concept of nonlinear power-law scalings goes back to the Kolmogorov theory of fully developed turbulence in the sixties (see \cite{cf:Castaing,cf:Sch,cf:Sto,cf:Cas,cf:Fr} and references therein), introduced to render the intermittency effects in turbulence. Mandelbrot \cite{mandelbrot} came up with the first mathematical discrete approach of multifractality, the now celebrated multiplicative cascades. Roughly speaking, a (dyadic) multiplicative cascade is a positive random measure $M$ on the unit interval $[0,1]$ that obeys the following decomposition rule:
\begin{equation}\label{star:casc}
M(dt)\stackrel{law}{=}Z^0\ind_{[0,\frac{1}{2}]}(t)M^0(2dt)+Z^1\ind_{[\frac{1}{2},1]}(t)M^1(2dt-1),
\end{equation}
where $M^0,M^1$ are two independent copies of $M$ and $(Z^0,Z^1)$ is a random vector with prescribed law and positive components of mean $1$ independent from $M^0,M^1$. Such an equation (and its generalizations to $b$-adic trees for $b\geq 2$), the celebrated star equation introduced by Mandelbrot in \cite{mandelbrotstar}, uniquely determines the law of the multiplicative cascade.
Despite the fact that multiplicative cascades have been widely used as reference models in many
applications, they possess many drawbacks related to their discrete scale invariance, mainly they
involve a particular scale ratio and they do not possess stationary fluctuations (this
comes from the fact that they are constructed on a dyadic tree structure).

Much effort has been made to develop a continuous parameter theory of suitable stationary multifractal random measures ever since, stemming from the theory of multiplicative chaos introduced by Kahane \cite{cf:Kah,Bar,cf:Sch,bacry,cf:RoVa,rhovar}. The construction of such measures is now well understood and they are largely used in mathematical modeling since they obey a so-called stochastic scale invariance property, namely the property of being equal in law at different scales up to an independent stochastic factor. This is some kind of continuous parameter generalization of (\ref{star:casc}).
Stochastic scale invariance property is observed in many experimental and theoretical problems, like turbulence (see \cite{cf:Fr,cf:Castaing}  and many others), quantum gravity (see \cite{cf: KPZ,cf:DuSh,cf:RhoVar}), mathematical finance, etc... and this is the main motivation for introducing multifractal random measures. However, as far as we know, the following question has never been solved: are these measures the only existing stochastic scale invariant object? This is fundamental since a positive answer gives a full justification to their intensive use. In this paper, we characterize stochastic scale invariant measures when the stochastic factor is assumed to be log-normal. We prove that the class of such objects is made up of Gaussian multiplicative chaos with a specific structure of the covariation kernel, which turns out to be larger than described in the literature.

\section{Background}
 Let us first remind the reader of the main definitions we will use throughout the paper. We denote by $\mathcal{B}(E)$ the Borelian sigma field on a topological space $E$.
 A random  measure $M$  is a random variable taking values into the set of positive Radon measures defined on $\mathcal{B}(\R^d)$ such that $\E[M(K)]<+\infty$ for every compact set $K$. 
A random measure $M$ is said to be stationary if for all $y\in\R^d$ the random measures $M(\cdot)$ and $M(y+\cdot) $ have the same law. 
\subsection{Gaussian multiplicative chaos}
We remind the reader of the notion of Gaussian multiplicative chaos as introduced by Kahane \cite{cf:Kah}. Consider a sequence $(X^n)_n$ of independent centered stationary Gaussian processes with associated nonnegative covariance kernel $k_n(r)=\E[X^n_rX^n_0] \geq 0$. For each $N\geq 1$, we can define a Radon measure $M^N$ on the Borelian subsets of $\R^d$ by
$$M^N(A)=\int_Ae^{\sum_{n=0}^NX^n_r-\frac{1}{2}\E[(X^n_r)^2]}\,dr.$$
For each Borelian set $A$, the sequence $(M^N(A))_N$ is a positive martingale. Thus it converges almost surely towards a random variable denoted by $M(A)$. One can deduce that the sequence of measures $(M^N)_N$ weakly converges towards a Radon measure $M$, commonly denoted by 
\begin{equation}\label{GMC}
M(A)=\int_Ae^{X_r-\frac{1}{2}\E[X_r^2]}\,dr 
\end{equation}
and called Gaussian multiplicative chaos associated to the kernel 
\begin{equation}\label{decompo}
K(r)=\sum_{n=0}^{+\infty}k_n(r).
\end{equation} 
Roughly speaking, (\ref{GMC}) can be understood as a measure admitting as density the exponential of a Gaussian process $X$ with covariance kernel $K$. Of course, this is purely formal because $X$ can only be understood as a (random Gaussian) distribution in the sense of Schwartz because of the possible singularities of the kernel $K$. 

Of special interest is the situation when the function $K$ can be rewritten as (for some $\lambda^2>0$)
\begin{equation}\label{equationexplicite}
K(r)=\lambda^2\ln_+\frac{T}{|r|}+g(r)
\end{equation}
for some bounded function $g$ (and $\ln_+(x)=\max(0,\ln(x))$). In that case, Kahane proved that the martingale $(M^N(A))_N$, for some Borelian set $A$ with non-null finite Lebesgue measure, is uniformly integrable if and only if $\lambda^2<2d$. This condition is necessary and sufficient in order for the limiting measure $M$ to be non identically null. For kernels of the form (\ref{equationexplicite}) which can not be written as a sum of nonnegative terms as (\ref{decompo}), we refer to the extended Gaussian multiplicative theory developed in \cite{cf:RoVa}. We remind that Gaussian multiplicative chaos with kernel given by (\ref{equationexplicite}) has found applications in many fields in science:
\begin{itemize}
 \item
 In dimension $1$, the measure $M$ is called the lognormal Multifractal Random Measure (MRM). It is used to model the volatility of a financial asset (see \cite{cf:BaKoMu}, \cite{cf:DuRoVa}). 
 \item
 In dimension $2$, the measure $M$ is a probabilistic formulation of the quantum gravity measure (more precisely, the quantum gravity measure is defined as the exponential of the Gaussian Free Field and therefore is defined in a bounded domain). We refer to references  \cite{Benj}, \cite{cf:DuSh}, \cite{cf:RhoVar} for probabilistic papers on this topic.  
 \item   
 In dimension $3$, the measure $M$ is called the Kolmogorov-Obhukov model (see textbook \cite{cf:Fr}): it is a model of energy dissipation in the statistical theory of fully developed turbulence.   
    
\end{itemize}

\section{Main results}

\subsection{Definitions}

In this paper we are interested in stationary random measures  satisfying the following scale invariance property:

\begin{definition}{\bf Log-normal $\star$-scale invariance.} \label{def_1}
A random measure $M$ is said to be lognormal $\star$-scale invariant if for all $\epsilon<1$, $M$ obeys the cascading rule
\begin{equation}\label{star}
 \big(M(A)\big)_{A\in\mathcal{B}(\R^d)}\stackrel{law}{=} \big(\int_Ae^{\omega_{\epsilon}(r)}M^{\epsilon}(dr)\big)_{A\in\mathcal{B}(\R^d)}
\end{equation}
where $\omega_{\epsilon}$ is a stationary Gaussian process with continuous sample paths and $M^{\epsilon}$ is a random measure independent from $\omega_{\epsilon}$ satisfying the relation
\begin{equation}\label{starbis}
\big(M^{\epsilon}(\epsilon A)\big)_{A\in\mathcal{B}(\R^d)}\stackrel{law}{=}\epsilon^d \big(M(A)\big)_{A\in\mathcal{B}(\R^d)}.
\end{equation}\qed
\end{definition}
Intuitively, this relation means that when you zoom in the measure $M$, you should observe the same behavior up to an independent log-normal factor. This relation is the continuous parameter analog of the celebrated Mandelbrot star equation.
\begin{remark}
In order for a measure $M$ satisfying (\ref{star}) with a moment of order $1$ to be non trivial, it is obvious to check that the Gaussian process $\omega_{\epsilon}$ must be normalized so that $\E[e^{\omega_{\epsilon}}(r)]=1$.
\end{remark}

\begin{definition}\label{goodness}
We will say that a stationary random measure $M$ satisfies the good lognormal $\star$-scale invariance  if $M$ is lognormal $\star$-scale invariant and  for each $\epsilon<1$, the covariance kernel $k_{\epsilon}$ of the 
process $\omega_{\epsilon}$ involved in (\ref{star})  is continuous and satisfies:
\begin{align}
 & &|k_{\epsilon}(r)|&\to 0\quad \text{ as }\quad |r|\to+\infty,\label{lim}\\
\forall r,r'\in\R^d\setminus\{0\}, & & |k_{\epsilon}(r)-k_{\epsilon}(r')|&\leq  C_{\epsilon}\theta\big(\min(|r|,|r'|)\big)|r-r'|\label{lipschitz}
\end{align} 
 for some positive constant $C_{\epsilon}$ and some  decreasing function $\theta:]0,+\infty[\to \R_+$ such that 
\begin{equation}\label{cvint}
\int_1^{+\infty}\theta(u)\ln(u)\,du<+\infty.
\end{equation}
\qed
\end{definition}

Though we would like to solve \eqref{star} in great generality, we must make a few technical assumptions to avoid pathological situations (a pathological example is given at the very end of Section \ref{built}). This is basically the purpose of the above definition \ref{goodness}.  Let us make a few comments on its content. 

Equation \eqref{lipschitz} mainly expresses that the kernel $k_{\epsilon}$ is Lipschitzian with a local Lipschitz constant that decays at most like $\theta$ when approaching infinity. 
By combining \eqref{lim} and \eqref{lipschitz}, it is plain to see that
\begin{align}
\forall r\not=0, & & |k_{\epsilon}(r)|&\leq C_{\epsilon}\int_{|r|}^{+\infty}\theta(u)\,du.\label{modulus}
\end{align} This is a very weak decorrelation property for the process $\omega_\epsilon$, which describes how fast the covariance function decays at infinity. In our proofs, it will be the key tool to investigate the mixing properties of the measure $M$.

\subsection{Results}

In what follows, we are mainly interested in the one-dimensional case  $d=1$. We  have the following description of the solutions to \eqref{star}, which is the main result of the paper:
\begin{theorem}\label{main}
Let $M$ be a   good lognormal $\star$-scale invariant random measure. Assume that $$\E[M([0,1])^{1+\delta}]<+\infty$$ for some $\delta>0$. Then $M$ is the product of a nonnegative random variable $Y\in L^{1+\delta}$  and an independent Gaussian multiplicative chaos 
\begin{equation}\label{chaos}
\forall A\subset \mathcal{B}(\R),\quad M(A)=Y\int_Ae^{X_r-\frac{1}{2}\E[X_r^2]}\,dr
\end{equation}
 with associated  covariance kernel given by the improper  integral
\begin{equation}\label{struct}
K(r)=\int_{|r|}^{+\infty}\frac{k(u)}{u}\,du
\end{equation} for some continuous covariance function $k$ such that   $k(0)\leq \frac{2}{1+\delta}$. 

Conversely, given some datas $k$ and $Y$ as above, the relation \eqref{chaos} defines a  log-normal $\star$-scale invariant random measure $M$ with finite  moments of order $1+\gamma$    for every $\gamma \in [0, \delta)$. 
\end{theorem}

Let us also state the following result giving a sufficient  (and not far from being necessary) condition in terms of $k$ for the measure $M$ as constructed in Theorem \ref{main} to be good:
\begin{proposition}\label{propgood}
Let  $M$ be a  log-normal $\star$-scale invariant random measure as constructed  in Theorem \ref{main}. If  
\begin{equation}\label{intln}
\int_1^{+\infty}\ln r\sup_{|u|\geq r}\frac{|k(u)|}{u}\,dr<+\infty 
\end{equation}  then $M$ is a good lognormal $\star$-scale invariant random measure.
\end{proposition}
 
Let us  comment on Theorem \ref{main}. First we point out that $Y$ is deterministic as soon as the random measure $M$  is ergodic. Second, good lognormal $\star$-scale invariant measures exhibit a multifractal behaviour. More precisely, if we consider a measure $M$ as in Theorem \ref{main}, we define its structure exponent
$$\forall q >0,\quad \xi(q)=(1+\frac{k(0)}{2})q-\frac{k(0)}{2}q^2.$$
Then we have the following asymptotic power-law spectrum, for $q< 1+\delta$:
$$\E\big[M([0,t])^q\big]\simeq C_q t^{\xi(q)}\quad \text{ as }t\to 0, $$ for some positive constant $C_q$.

We also stress that the intermittency parameter $k(0)$ is explicit when one knows $K$ because of the relation
\begin{equation}
K(r) \sim k(0) \ln\left(\frac{1}{r}\right), \quad \text{ when } r\to 0.
\end{equation}
The covariance function $K$ can also be recovered from the two sets marginals of the measure $M$   thanks to formula \eqref{limh}. 
 
Finally, Theorem \ref{main} has the following consequence about the regularity of good lognormal $\star$-scale invariant measures:
\begin{corollary}\label{atom}
Almost surely, a good log-normal $\star$-scale invariant random measure $M$ does not possess any atom on $\R$, that is:
$$\mathrm{almost}\,\, \mathrm{surely}, \forall x\in \R,\quad M(\{x\})=0.$$
\end{corollary}

 Now we investigate  long-range independence for good lognormal $\star$-scale invariant random measures.
 So we  introduce the related notion of cut-off:
\begin{definition}\label{def_cut-off}
We will say that a stationary random measure $M$ admits a cut-off $d>0$ if, for $t<s$, the $\sigma$-algebras $\mathcal{H}_{-\infty}^t=\sigma\{M(A);A\in \mathbb{B}(\R), A\subset (-\infty,t]\}$ and  $\mathcal{H}_{s}^{+\infty}=\sigma\{M(A);A\in \mathbb{B}(\R), A\subset [s,+\infty)\}$ are independent, conditionally to the asymptotic $\sigma$-algebra of $M$, as soon as $s-t>d$.\qed
\end{definition}
 Of course, if the measure $M$ is ergodic then  the asymptotic $\sigma$-algebra of $M$ is trivial and we can remove the sentence "conditionally to the asymptotic $\sigma$-algebra of $M$" from the definition. For instance the measure constructed in subsection \ref{example} admits a cut-off $T$ and is ergodic. It results from the proof of Theorem \ref{main} that the cut-off property can be read off the cascading rule \eqref{star}:
 
 \begin{proposition}\label{propcut}
Let $M$ be a good lognormal $\star$-scale invariant  random measure with finite $1+\delta$ moment. Then $M$ admits a cutoff if and only if,  for some $\epsilon <1$ (or equivalently for all $\epsilon <1$), the covariance kernel $k_{\epsilon}$ of the process $\omega_{\epsilon}$ in \eqref{star} reduces to $0$ outside a compact set.
 \end{proposition}
 
 Finally, we mention that another notion of stochastic scale invariance has been studied in the literature before: it is called the exact stochastic scale invariance 
 (see \cite{bacry,cf:Castaing,rhovar}). Let us recall the main result:  if the Gaussian multiplicative chaos $M$ admits a covariance kernel $K$ such that $K(x)=\lambda^2\ln \frac{T}{|x|}+C$ 
 for some constant $C$ and for all $x$ in a ball $B(0,R)$ then $M$ satisfies the "exact stochastic scale invariance":
 $$\forall \alpha\in (0,1),\quad (M(\alpha A))_{A\subset B(0,R)}\stackrel{law}{=}\alpha e^{Y_\alpha-\frac{1}{2}\E[Y_\alpha^2]} (M( A))_{A\subset B(0,R)} $$ where $Y_\alpha$ is a 
 centered Gaussian random variable with variance $\lambda^2\ln \frac{1}{\alpha} $. 
 
The reader may wonder if we can construct random measures that are both exactly stochastically scale invariant and good lognormal $\star$-scale invariant. Let us 
 show that 
 
 \begin{proposition}\label{exactsi}
 Let $M$ be a Gaussian multiplicative chaos whose covariance kernel $K$ is such that, for $|r|\leq R$, $K(r)=\lambda^2\ln \frac{T}{|r|}+C$ 
 for some constant $C$ (in particular, $M$ satisfies the "exact stochastic scale invariance"), then $M$ is not a good lognormal $\star$-scale invariant 
 random measure. 
 \end{proposition}

\subsection{Multidimensional results}
We stress that our results remain true in higher dimensions without changes in the proofs.  For the sake of completeness, we state the main result.  

\begin{theorem}\label{multi}
Let $M$ be a   good lognormal $\star$-scale invariant random measure such that for each $\epsilon<1$, the covariance kernel $k_{\epsilon}$ of the process $\omega_{\epsilon}$ is continuous and differentiable on $\R^d\setminus\{0\}$. Assume that $$\E[M([0,1]^d)^{1+\delta}]<+\infty$$ for some $\delta>0$. Then $M$ is the product of a nonnegative random variable $Y\in L^{1+\delta}$  and an independent Gaussian multiplicative chaos: 
\begin{equation}\label{chaosmulti}
\forall A\subset \mathcal{B}(\R^d),\quad M(A)=Y\int_Ae^{X_r-\frac{1}{2}\E[X_r^2]}\,dr
\end{equation}
 with associated  covariance kernel given by the improper integral
\begin{equation}\label{structmulti}
\forall x\in \R^d\setminus\{0\},\quad K(x)=\int_{1}^{+\infty}\frac{k(x u)}{u}\,du
\end{equation} for some continuous covariance function $k$ such that   $k(0)\leq \frac{2d}{1+\delta}$. 

Conversely, given some datas $k$ and $Y$ as above, the relation \eqref{chaos} defines a  lognormal $\star$-scale invariant random measure $M$ with finite  moments of order $1+\gamma$ for every $\gamma \in [0, \delta)$.  
\end{theorem}
 
 It turns out that Proposition \ref{propgood} remains true in dimension $d\geq1$. When the dimension is greater than $1$, it may be interesting to focus on the isotropy properties. In the same spirit as Proposition \ref{propcut}, for a good lognormal $\star$-scale invariant 
 measure $M$ with a finite moment of order $1+\delta$, the following assertions are equivalent:
\begin{enumerate}
\item $M$ is isotropic,
\item its covariance kernel $K$ (or equivalently $k$ in \eqref{structmulti}) is isotropic,
\item the covariance kernel $k_{\epsilon}$ is isotropic for some $\epsilon<1$,
\item the covariance kernels $k_{\epsilon}$ are isotropic for all $\epsilon <1$.
\end{enumerate}
\subsection{Classical example }\label{example}
As far as we know, there exists only one example of good log-normal $\star$-scale invariant random measures in the literature, which was first described in \cite{Bar} (see also \cite{bacry}). Its construction is very intuitive: it is geometric and relies on homothetic properties  of triangles in the half-plane. We also stress that this specific example of $\star$-scale invariant random measures is not restricted to  the Gaussian case:  the factor can be more general (log-L\'evy).  

Following \cite{bacry}, we recall the construction of this example and refer the reader to the aforementioned papers for further details. Fix $T>0$ and let $\mathcal{S}^+$ be the state-space half plane
\begin{equation*}
\mathcal{S}^+ = \{(t,l): t\in \R, l>0\}. 
\end{equation*}
with which one can associate the measure 
\begin{equation*}
\mu(dt,dl) = l^{-2} dt dl. 
\end{equation*}
Then we introduce the independently scattered Gaussian random measure $P$ defined for any $\mu$-measurable set $A$ by
\begin{equation*}
\E\left[e^{iq P(A)}\right] = e^{\varphi(q)\mu(A)} 
\end{equation*}
with $\varphi(q)=-\lambda^2 q^2 /2 - i q \lambda^2/ 2$. Under those assumptions, we can note that for any $\mu$-measurable set $A$, $P(A)$ is a Gaussian 
variable with mean $m=-\mu(A)\lambda^2/2$ and variance $\sigma^2=\lambda^2 \mu(A)$.
We can then define the Gaussian process $(\omega_l(t))_{t\in \R}$ for $l \geq 0$ by
\begin{equation*}
\omega_l(t) = P\left(\mathcal{A}_l(t)\right)
\end{equation*}
where $\mathcal{A}_l(t)$ is the triangle like subset $\mathcal{A}_l(t):=\{(t',l'): l \leq l' \leq T, -l'/2 \leq t-t' \leq l'/2\}$. 

\begin{figure}[h]
\begin{tikzpicture}[domain=0:7]
\draw[->] (-0.1,0) -- (8.1,0) node[pos=0,below]{0} node[pos=0.5,below]{$t$} ;
\draw[->] (0,-0.2) -- (0,5.2)   node[above] {$l$};
 \draw[color=gray,-,thick, fill=gray] (1,4) -- (7,4) node[right,thick,black]{$\mathcal{A}_l(t)$} -- (4.75,1) -- (3.25,1) -- (1,4);
 \draw[color=gray,style=dashed] (3.25,1) -- (4,0) -- (4.75,1);
  \draw[color=gray,style=dashed] (3.25,1) -- (0,1) node[left,color=black]{$l$};
   \draw[color=gray,style=dashed] (1,4) -- (0,4) node[left,color=black]{$T$};
   
\end{tikzpicture}
\end{figure}

Define now the random measure $M_l$ by $M_l(dt) = e^{\omega_l(t)} dt$. Almost surely, the family of measures $(M_l(dt))_{l>0}$ weakly converges towards a  random measure $M$. If $\lambda^2<2$, this measure is not trivial. 

Let us check that $M$ is  a good log-normal $\star$-scale invariant random measure. Fix $\epsilon<1$
and define the sets $\mathcal{A}_{l,\epsilon T}(t):=\{(t',l'): l \leq l' \leq \epsilon T, -l'/2 \leq t-t' \leq l'/2\}$ and 
$\mathcal{A}_{\epsilon T, T}(t):=\{(t',l'): \epsilon T \leq l' \leq T, -l'/2 \leq t-t' \leq l'/2\}$. Note that $\mathcal{A}_l(t) = \mathcal{A}_{l,\epsilon T}(t) \cup \mathcal{A}_{\epsilon T,T}(t) $
and that those two sets are disjoint. 
Thus, we can write for every $\mu$-measurable set $A$
\begin{equation}\label{eq_l}
M_l(A) = \int_A e^{\omega_{\epsilon T,T}(t)} e^{\omega_{l,\epsilon T}(t)} dt 
\end{equation}
with $\omega_{\epsilon T,T}(t) = P(\mathcal{A}_{\epsilon T,T}(t))$ and $\omega_{l,\epsilon T}(t)= P( \mathcal{A}_{l,\epsilon T}(t))$.

\begin{figure}[h]
\begin{tikzpicture}[domain=0:7]
\draw[->] (-0.1,0) -- (8.1,0) node[pos=0,below]{0} node[pos=0.5,below]{$t$} ;
\draw[->] (0,-0.2) -- (0,5.2)   node[above] {$l$};
 \draw[color=blue!60,-,thick, fill=blue!60] (1,4) -- (7,4) node[right,thick,blue]{$\mathcal{A}_{\epsilon T,T}(t)$} -- (5.5,2) -- (2.5,2) -- (1,4);
 \draw[color=red!60,-,thick, fill=red!60] (5.5,2) -- (2.5,2)  -- (3.25,1) -- (4.75,1) node[right,thick,red]{$\mathcal{A}_{l,\epsilon T}(t)$} -- (5.5,2);
 \draw[color=gray,thin,style=dashed] (3.25,1) -- (4,0) -- (4.75,1);
  \draw[color=gray,thin,style=dashed] (2.5,2) -- (0,2) node[left,color=black]{$\epsilon T$};
  \draw[color=gray,thin,style=dashed] (3.25,1) -- (0,1) node[left,color=black]{$l$};
   \draw[color=gray,thin,style=dashed] (1,4) -- (0,4) node[left,color=black]{$T$};
   
\end{tikzpicture}
\end{figure}

We then study equation (\ref{eq_l}) in the limit $l \to 0$; we obtain
\begin{equation}
M(A) = \int_A e^{\omega_{\epsilon T,T}(t)}  M^{\epsilon}(dt)
\end{equation}
where $M^{\epsilon}$ is the limit when $l \to 0$ of the random measure $M^\epsilon_l(dt):=e^{\omega_{l,\epsilon T}(t)} dt$. 
We easily verify that $M^\epsilon(\epsilon A) \stackrel{law}{=} \epsilon M(A)$ writing
\begin{equation}
M^\epsilon_l(A) = \epsilon \int_A e^{\omega_{l,\epsilon T}(\epsilon t)}dt 
\end{equation}
and checking that the covariance of the Gaussian process $(\omega_{l,\epsilon T}(\epsilon t))_{t \in \R}$ is the same as the one of $(\omega_{l,T}(t))_{t\in \R}$. 

The covariance kernel of the stationary Gaussian process $\omega_{\epsilon T,T}(t)$ is given by
\begin{equation}
k_\epsilon(r) = 
\left\{ 
\begin{array}{ll}
0 & \qquad \mathrm{if} \quad |r| \geq T \\ 
\lambda^2 (\ln \frac{T}{|r|} + \frac{|r|}{T} - 1 ) & \qquad \mathrm{if} \quad \epsilon T \leq |r| \leq T \\
\lambda^2 (\ln \frac{1}{\epsilon} + \frac{|r|}{T} - \frac{|r|}{\epsilon T}) &\qquad \mathrm{if} \quad |r|\leq \epsilon T\,. 
\end{array}
\right.
\end{equation}
Since $k_\epsilon$ reduces to $0$ outside a compact set, it is straightforward to check \eqref{lim} and \eqref{lipschitz}. We further stress that this measure admits a cut-off in the sense 
of Definition \ref{def_cut-off}.  

\begin{remark}
In view of Theorem \ref{main}, note that the random measure $M$ is a Gaussian multiplicative chaos with associated kernel 
\begin{equation}\label{barral}
K(r)=\int_{|r|}^{+\infty}\frac{k(u)}{u}\,du \quad \text{with }\quad k(u)=\lambda^2(1-\frac{|u|}{T})\ind_{[0,T]}(|u|).
\end{equation}
and that we have 
$$k_{\epsilon}(r)=\int_{|r|}^{\frac{|r|}{\epsilon}}\frac{k(u)}{u}\,du.$$
\end{remark}

\section{Construction of log-normal $\star$-scale invariant random measures}\label{built}
This section is devoted to the existence part of Theorem \ref{main}: we give an explicit construction of lognormal $\star$-scale invariant random measures. 

We are given a positive random variable $Y\in L^{1+\delta}$ (for some $\delta>0$) and
a continuous covariation kernel $k$ such that $k(0)\leq \frac{2}{1+\delta}$. Let $F$ be the  (symmetric) spectral measure associated to $k$, that is
$$k(t)=\int_\R e^{i\lambda t}F(d\lambda),$$ and we assume that the improper integral $$K(r)=\int_r^{+\infty}\frac{k(u)}{u}\,du$$ converges for $r>0$. 

Let $\mu,\nu$ be two i.i.d. independently scattered Gaussian random measures (independent of $Y$) distributed on the half plane $\R\times \R^*_+$ such that:
$$\forall A\in \mathcal{B}(\R\times \R^*_+),\quad \E[e^{q\mu(A)}]=e^{\frac{1}{2}q^2\theta(A)}$$
where  
$$\theta(A)=\int_{\lambda\in \R}\int_{y\in \R_+^*}\ind_A(\lambda,y)\frac{1}{y}dy F(d\lambda).$$
Let $\epsilon<1$, we define the centered Gaussian process
$$\forall t \in \R, \quad X_{\epsilon}(t)=\int_{\lambda\in \R}\int_{y\in [1,\frac{1}{\epsilon}[}\cos(\lambda ty )\mu(d\lambda,dy)+\int_{\lambda\in \R}\int_{y\in [1,\frac{1}{\epsilon}[}\sin(\lambda ty )
\nu(d\lambda,dy).$$
It is plain to compute its covariation kernel, call it $k_{\epsilon}$, by using the symmetry of the spectral measure $F(d\lambda)$:
\begin{align*}
k_{\epsilon}(t-s) &=\E[X_{\epsilon}(s)X_{\epsilon}(t)]\\ &=\int_{\lambda\in \R}\int_{y\in [1,\frac{1}{\epsilon}[}\cos(\lambda ty )\cos(\lambda t s)\frac{1}{y}dy F(d\lambda)+
\int_{\lambda\in \R}\int_{y\in [1,\frac{1}{\epsilon}[}\sin(\lambda ty )\sin(\lambda sy )\frac{1}{y}dy F(d\lambda)\\
&=\int_{\lambda\in \R}\int_{y\in [1,\frac{1}{\epsilon}[}\cos(\lambda (t-s)y )\frac{1}{y}dy F(d\lambda)\\
&=\int_{y\in [1,\frac{1}{\epsilon}[}\int_{\lambda\in \R}e^{i\lambda (t-s)y }F(d\lambda)\frac{1}{y}dy \\
&=\int_{y\in [1,\frac{1}{\epsilon}[}\frac{k(|t-s|y)}{y}dy \\
&=\int_{|t-s|}^{\frac{1}{\epsilon}|t-s|}\frac{k(y)}{y}dy.
\end{align*}

For all $A\in \mathcal{B}(\R)$, the process $$M_{1/l}(A)=Y\int_{A} \exp \big(X_{{1/l}}(r)-\frac{1}{2}\E[X^2_{{1/l}}(r)]\big) dr$$ is  obviously a positive martingale and thus 
converges as $l \to \infty$ towards a random variable $M(A)$. The stationary random measure $(M(A))_{A\in \mathcal{B}(\R)}$ is a Gaussian multiplicative chaos in the 
sense of \cite{cf:RoVa} with associated kernel $K$.

Note that for $l>1/\epsilon$, we have $\forall t \in \R$:
\begin{align}
 X_{1/l}(t)=& X_{\epsilon}(t)+\int_{\lambda\in \R}\int_{y\in [\frac{1}{\epsilon},l[}\cos(\lambda ty )\mu(d\lambda,dy)+
 \int_{\lambda\in \R}\int_{y\in [\frac{1}{\epsilon},l[}\sin(\lambda ty )\nu(d\lambda,dy)\nonumber\\
 \stackrel{def}{=}& X_{\epsilon}(t)+\bar{X}_{\epsilon,1/l}(t),\label{xl}
\end{align}
where $\bar{X}_{\epsilon,1/l}$ is a centered stationary Gaussian process independent from $X_{\epsilon}$ with covariance kernel given by:
$$\bar{k}_{\epsilon,1/l}(t-s)= \E[\bar{X}_{\epsilon,1/l}(s)\bar{X}_{\epsilon,1/l}(t)]=\int^{l|t-s|}_{\frac{1}{\epsilon}|t-s|}\frac{k(y)}{y}dy .$$
As above, we can define the random measure $M^{\epsilon}$ as the limit as $l\to +\infty$ of the random measures
$$\forall A\in \mathcal{B}(\R),\quad M^{\epsilon}_{1/l}(A)=Y\int_{A} \exp \big(\bar{X}_{\epsilon,1/l}(r)-\frac{1}{2}\E[\bar{X}_{\epsilon,1/l}^2(r)]\big) dr.$$
The stationary random measure $(M^{\epsilon}(A))_{A\in \mathcal{B}(\R)}$ is a Gaussian multiplicative chaos in the sense of \cite{cf:RoVa} with associated 
covariance $K(\cdot \frac{1}{\epsilon})$. We deduce that $\frac{1}{\epsilon}M^{\epsilon}\big(\epsilon\cdot\big)$ is a  Gaussian multiplicative chaos in the sense 
of \cite{cf:RoVa} with associated covariance $K(\cdot)$. The measure $\frac{1}{\epsilon}M^{\epsilon}\big(\epsilon\cdot\big)$ thus has the same law as $M$. From \eqref{xl}, we obviously have:
$$M(A)=\int_A\exp \big(X_{\epsilon}(r)-\frac{1}{2}\E[X_{\epsilon}^2(r)]\big) M^{\epsilon}(dr)$$ in such a way that \eqref{star} holds.
Finally we point out that   $M$ admits a moment of order $1+\gamma$ for all $0\leq \gamma<\delta$ (see \cite{cf:Kah}).

\begin{remark}
By focusing on the above construction, we see that the covariance kernel $k$ can be intuitively interpreted as some kind of infinitesimal stochastic generator. We may look
 $X_\epsilon$ as a sum
 $$X_\epsilon(r)=\sum_{1\leq y \leq \frac{1}{\epsilon}}a_yZ^y_r$$ where $(Z^y)_y$  are independent centered Gaussian processes with kernel $k(y\cdot)$ and $(a_y)_y$ are independent random Gaussian variables with variance $\frac{dy}{y}$. So, when $\epsilon$ decreases infinitesimally, we "add" an independent Gaussian process with kernel $k(\frac{1}{\epsilon}\cdot)$ times an independent Gaussian factor of variance $\frac{-d\epsilon}{\epsilon}$.
\end{remark}

\noindent {\it Proof of Proposition \ref{propgood}}
We show that the measure $M$ is good under assumption \eqref{intln}. Because $k$ is continuous, the kernel $k_\epsilon(r)=\int_{|r|}^{|r|/\epsilon}\frac{k(u)}{u}\,du$ is of class $C^1$ on $\R^*$. Thus, we have:
$$|k_\epsilon(r)-k_\epsilon(r')|\leq \sup_{u\geq \min(|r|,|r'|)}|k_\epsilon'(u)|.$$
Because we have $$k_\epsilon'(u)=\frac{1}{u}(k(u/\epsilon)-k(u)),$$
it is plain to see that a reasonable choice for $\theta$ is $\theta(x)=\sup_{u\geq|x|}|\frac{k(u)}{u}|$ and $C_\epsilon=2/\epsilon$.
$$\int_1^{+\infty}\ln r\sup_{|u|\geq r}\frac{|k(u)|}{u}\,dr<+\infty\Rightarrow\int_1^{+\infty}\ln r\,\theta(r)\,dr<+\infty,$$ so that the measure is good. \qed

\subsection{Practical examples}
In this subsection, we give practical examples of log-normal $\star$-scale invariant random measures. Using Theorem \ref{main}, 
{\bf good} log-normal  $\star$-scale invariant random measures are Gaussian multiplicative chaos whose covariance structure is given by
\begin{equation}
K(s) = \int_{|s| }^{+\infty} \frac{k(u)}{u} du 
\end{equation}
where $k$ is a continuous covariance function satisfying $k(0)<2$ and some weak decay assumptions (ensuring \eqref{intln} for instance). Therefore, to define explicit examples, we just need to exhibit suitable kernels $k$. The decay assumptions can be read off the spectral measure of $k$. For instance, if $k$ is the Fourier transform of some positive even integrable function $f$, which possesses an integrable derivative, it is a simple application of the Riemann theorem to prove that \eqref{intln} is satisfied. Actually, for \eqref{intln} to be satisfied, the assumptions on the regularity of the spectral measure can be much weakened. For instance, we can consider a kernel $k$ that is the Fourier transform of some positive even integrable function $f$ with integrable $\alpha$-fractional derivative for $0<\alpha<1$:
$$\partial_\alpha f=\int_{\R^*}\frac{f(x+z)-f(x)}{|z|^{1+\alpha}}\,dz\in L^1(\R).$$
In that case, the Riemann theorem implies $|u|^\alpha k(u)\to 0$ as $|u|\to \infty$ and it is then plain to se that \eqref{intln} is satisfied.

Below are listed a few examples of such kernels: 
\begin{itemize}
\item the function $k(s)  = \frac{1}{\sigma\sqrt{2\pi}} e^{-\frac{|s|^2}{2\sigma^2}}$ (where $\sigma>0$) is continuous and positive-definite since its Fourier transform 
$\hat{k}(r)= e^{-\sigma^2 r^2/2}$ is positive. 
\item the covariance function of the stationary Orstein-Uhlenbeck process which takes on the form $k(s) = \frac{\sigma^2}{2\theta} e^{-\theta |s|}$ where $\theta>0,\sigma>0$.
\item we can consider $k$ as the Fourier transform of the function ( $\lambda>0$)
\begin{equation}\label{resolvent}
f(x)=\int_0^{+\infty}e^{-\lambda t}\E[g(x+X_t)]\,dt
\end{equation}
 where $g\in L^1(\R)$ is any positive integrable function and $X$ is a pure jump L\'evy process with L\'evy symbol 
$$\eta(u)=\int_{\R^*}(e^{iuz}-1)\frac{1}{|z|^{1+\alpha}}\,dz $$ for some $0<\alpha<1$. It is well know that the Lebesgue measure is invariant for the semi-group generated by $X$ so that $k(0)=\|f\|_1=\|g\|_1/\lambda$: this gives a condition on the norm $\|g\|_1$ for having $k(0)<2$. Furthermore, $f$ admits an integrable $\alpha$-fractional derivative so that \eqref{intln} is satisfied. Actually, it turns out that all the functions in $L^1(\R)$ with an integrable $\alpha$-fractional derivative admit a representation as \eqref{resolvent}. The reader may consult \cite{applebaum} for further details.
\end{itemize}
 
We stress that, as soon as they are not trivial (i.e. $k(0)<2$), the Gaussian multiplicative chaos of the first two above examples do not have cut off in the sense of Definition \ref{def_cut-off}.  Obviously, many other examples exist. 

Let us mention another example of log-normal $\star$-scale invariant random measures which does not present the goodness property of Definition \ref{goodness}.  
From Theorem \ref{main}, the Gaussian multiplicative chaos associated to the covariance function
\begin{equation}
K(s)  = \int_{|s|}^{+\infty} \frac{\cos(u)}{u} du. 
\end{equation}
is log-normal $\star$-scale invariant  in the sense of Definition \ref{def_1}. 
The function $k(r) = \cos(r)$ is indeed positive definite since its spectral measure is the positive measure $(\delta_{1}(dx)+\delta_{-1}(dx))/2$. The kernel $k_\epsilon(r)=\int_{|r|}^{|r|/\epsilon}\frac{\cos(u)}{u} du$ does not satisfy \eqref{lipschitz} so that the associated measure $M$ is not good.
Note that  this Gaussian multiplicative chaos falls under the scope of \cite{cf:RoVa} since   the function $K$ does not have a constant positive sign. 


\section{Characterization of the good log-normal $\star$-scale invariant random measures}
This section is devoted to the proof of the first statement of Theorem \ref{main}. For the sake of readability, some proofs of auxiliary results are gathered in the appendix.
 
Let $M$ be a good log-normal scale
invariant random measure  defined on a probability space $(\Omega,\mathcal{F},\p)$. We introduce as usually the spaces $L^p$ 
on $(\Omega,\mathcal{F},\p)$ for $1\leq p\leq \infty$. Recall that the measure $M$ satisfies, for all $\epsilon\in (0,1)$
\begin{equation}\label{star2}
 \big(M(A)\big)_{A\in\mathcal{B}(\R)}\stackrel{law}{=} \big(\int_Ae^{\omega_{\epsilon}(r)}M^{\epsilon}(dr)\big)_{A\in\mathcal{B}(\R)}
\end{equation}
where $\omega_{\epsilon}$ is a Gaussian process independent from $M^{\epsilon}$, with $M^\epsilon(dr)=\epsilon M(\frac{dr }{\epsilon})$ in law. $k_\epsilon$ denotes the covariation kernel 
of the process $\omega_\epsilon$. Furthermore, we assume that the measure $M$ is non trivial ($M\not=0$) with a moment of order $1+\delta$ so that the process $\omega_\epsilon$ is necessarily normalized, 
that is $\E[e^{\omega_\epsilon}]=1$.

Now we introduce some definitions and tools that will be used throughout this section.
For each $\epsilon\in (0,1)$, define
\begin{equation}\label{ke}
\forall r\not=0,\quad K^\epsilon(r)=\sum_{n=0}^{+\infty}k_\epsilon\big(\frac{r}{\epsilon^n}\big).
\end{equation}
The uniform convergence of the series  on the sets $\{r\in\R;|r|>\rho\}$ for any $\rho>0$ is ensured by \eqref{modulus} since for $|r|>\rho$:
\begin{align}
\sum_{n=0}^{+\infty}|k_\epsilon\big(\frac{r}{\epsilon^n}\big)|\leq & C_\epsilon\sum_{n=0}^{+\infty}\int_\frac{|r|}{\epsilon^n}^{+\infty}\theta(u)\,du
\leq  C_\epsilon\sum_{n=0}^{+\infty}\int_\frac{\rho}{\epsilon^n}^{+\infty}\theta(u)\,du\nonumber\\
\leq & C_\epsilon\int_0^{+\infty}\int_{\rho\epsilon^{-y+1}}^{+\infty}\theta(u)\,du\,dy\nonumber\\
=& C_\epsilon\int_{\rho\epsilon}^{+\infty}\theta(u)\int_0^{\frac{\ln\frac{u}{\rho}}{-\ln \epsilon}+1}\,dy \,du\nonumber\\
 =& \frac{C_\epsilon}{-\ln \epsilon}\int_{\epsilon\rho}^{+\infty}\theta(u)\ln \frac{u}{\epsilon\rho}\,du
 \label{compk}
\end{align}
and this last integral is assumed to be converging \eqref{cvint}.  Furthermore, \eqref{lipschitz} also ensures that $K^\epsilon$ is Lipschitzian over each set $\{z\in\R;|z|>\rho\}$ for any $\rho>0$ because:
 \begin{align*}
|K^\epsilon(r)-K^\epsilon(r')|\leq &\sum_{n=0}^{+\infty}|k_\epsilon\big(\frac{r}{\epsilon^n}\big)-k_\epsilon\big(\frac{r'}{\epsilon^n}\big)|\\
\leq &C_\epsilon \sum_{n=0}^{+\infty}\theta\big(\frac{\min(|r|,|r'|)}{\epsilon^n}\big)\big|\frac{r-r'}{\epsilon^n}\big|\\
\leq & C_\epsilon\int_{0}^{+\infty}\theta\big(\frac{ \rho}{\epsilon^{y-1}}\big)\big|\frac{r-r'}{\epsilon^{y}}\big|\,dy\\
\leq & \frac{C_\epsilon}{-\rho \epsilon\ln \epsilon}|r-r'|\int_{\rho\epsilon}^{+\infty}\theta(u)\,du .
\end{align*}

We let $(X^n)_n$ denote a sequence of independent centered stationary Gaussian processes  with respective covariance kernels  $$\E[X^n_rX^n_s]=k_{\epsilon}(\frac{r-s}{\epsilon^{n}})\,\,\,\stackrel{\text{def}}{=}\overline{k}_n(r-s).$$ Clearly $X^n$ depends on $\epsilon$ but this parameter is omitted  from the notations for the sake of readability. We assume that the whole sequence $(X^n)_n$ and the measure $M$ are constructed on the same probability space and are mutually independent. We further define the measure $M^N$ for $N\geq 0$ by 
$$\forall A\in \mathcal{B}(\R),\quad M^N(A)=\epsilon^{N+1}M\big(\frac{1}{\epsilon^{N+1}}A\big).$$ 
Note that $\E[M^N(A)] = |A|$ where $|A|$ stands for the Lebesgue measure of the set $A$.

By iterating the scale invariance relation \eqref{star}, it is plain to see that, for each $N\geq 0$, the measure $\widetilde{M}^N$ defined by
\begin{equation}\label{mn}
\widetilde{M}^N(A)=\int_A\exp\Big(\sum_{n=0}^NX^n_r-\frac{1}{2}\E[(X^n_r)^2]\Big)\,M^N(dr)
\end{equation}
 has the same law as the measure $M$.

\subsection{Ergodic properties}
 First we investigate the immediate properties of $M$ resulting from the definitions.
\begin{lemma}\label{birk}
Let $M$ be a stationary random measure on $\R$ admitting a moment of order $1+\delta$. There is a nonnegative integrable random variable $Y\in L^{1+\delta}$ such that, for every bounded interval $I\subset \R$, $$\lim_{T \to \infty} \frac{1}{T} M\left(T I\right) = Y |I|\quad \text{almost surely and in }L^{1+\delta},$$
where  $|\cdot|$ stands for the Lebesgue measure on $\R$. As a consequence, almost surely the random measure $$A\in \mathcal{B}(\R)\mapsto \frac{1}{T}M(TA)$$ weakly converges towards $Y|\cdot|$ and $\E_Y[M(A)]=Y |A|$ ($\E_Y[\cdot]$ denotes the conditional expectation with respect to $Y$).
\end{lemma} 
 
 \noindent{\it{Proof.}}
If $M$ is a stationary random measure, the Birkhoff ergodic theorem implies the following convergence, for $n \in \N, n \to \infty$,
\begin{equation}
\frac{1}{n} M([0,n]) = \frac{1}{n} \sum_{i=1}^{n} M([i-1,i] )\to Y \quad \text{almost surely and in }L^{1+\delta}
\end{equation}
where $Y\in L^{1+\delta}$ is a nonnegative random variable. Using monotonicity of the mapping $t\mapsto M([0,t])$, one can show that $\frac{1}{T} M([0,T]) \to Y$ almost surely and in $L^{1+\delta}$.
For $a>0, b>a$, it is clear that  $\frac{1}{T} M\left(T[0,a]\right) \to a Y$ and that $\frac{1}{T} M\left(T[a,b]\right) \to (b-a) Y$ almost surely and in $L^{1+\delta}$.
So, for every bounded interval $I \subset \R_+$,  the following convergence holds  $\frac{1}{T} M(T I) \to |I| Y$ almost surely and in $L^{1+\delta}$.
Along the same lines, one can show the same convergence for every bounded interval $I \subset \R_-$ involving some nonnegative random variable $Y'\in L^{1+\delta}$.
Stationarity implies that $\frac{1}{T} M\left(T [-1,1]\right)$ has the same law as $\frac{1}{T} M\left(T[0,2]\right)$. By letting $T$ go to $\infty$, we find that $Y+Y'$ has the same law 
as $2Y$. Stationarity also implies that $Y'$ has the same law as $Y$.
 Let $0<\alpha<1$. We prove
\begin{equation}
\E[Y^\alpha] = \E\left[\left(\frac{Y+Y'}{2}\right)^\alpha\right] \\
                     \geq \frac{1}{2} \left( \E[Y^\alpha] + \E[Y'^\alpha] \right) \\
                    = \E[Y^\alpha]
\end{equation}
by using the Jensen inequality for the concave function $x\mapsto x^\alpha$. So the above inequality turns out to be  an equality and thus  $Y=Y'$ almost surely. We have shown that $\frac{1}{T} M(T I) \to |I| Y$ almost surely and in $L^{1+\delta}$ when $T \to \infty$  for every bounded interval $I \subset \R$.

Finally, by the portemanteau theorem, the convergence of the measure $A\in \mathcal{B}(\R)\mapsto \frac{1}{T}M(TA)$ on the intervals  towards $Y|\cdot|$ is enough to ensure the weak convergence.  \qed

\subsection{Mixing properties}

This section is devoted to study of the mixing properties of the measure $M$, which can be read off the structure of the kernel $K^\epsilon$.  

We first draw attention to the following relation, which will be used throughout the paper:
 $$\E_Y\big[F\big(M(A_1),\dots,M(A_n)\big)\big ]=\E_Y\big[F\big(\widetilde{M}^N(A_1),\dots,\widetilde{M}^N(A_n)\big)\big ] \quad a.s.$$
for every positive measurable function $F:\R^n\to \R$. The proof is deferred to appendix \ref{app:proof} (see Lemma \ref{lemY}).

\begin{lemma}\label{lemmix}
Let $A,B$ be two disjoint sets such that ${\rm dist}(A,B)>0$. Then the random variable $M(A)M(B)$ is integrable under $\E_Y[.]$ and
$$\E_Y[M(A)M(B) ] = Y^2\int_{A\times B}e^{K^\epsilon(r-u)}dr\,du.$$
\end{lemma}

\noindent {\it Proof.} We fix $R>0$ and denote by $\mathcal{G}$ the $\sigma$-field generated by $M$. Because the function $x\in \R_+\mapsto \min(R,x)$ is concave, we have
\begin{align*}
\E_Y\big[\min\big(R,M(A)M(B)\big)\big ] =&\E_Y\big[\min\big(R,\widetilde{M}^N(A)\widetilde{M}^N(B)\big)\big]\\
=&\E_Y\big[\E\big[\min\big(R,\widetilde{M}^N(A)\widetilde{M}^N(B)\big)|\mathcal{G}\big]\big]\\
\leq &\E_Y\Big[\min\Big(R,\E\big[\widetilde{M}^N(A)\widetilde{M}^N(B)|\mathcal{G}\big]\Big)\Big].
\end{align*}
Since $\widetilde{M}^N$ is given by \eqref{mn}, it is straightforward to compute:
\begin{equation}\label{ecs}
\E\big[\widetilde{M}^N(A)\widetilde{M}^N(B)|\mathcal{G}\big]=\int_{A\times B}e^{\sum_{n=0}^N\bar{k}_n(r-u)} M^N(dr)M^N(du).
\end{equation}
 Because of the uniform convergence of the series $\big(\sum_{n=0}^N\bar{k}_n(r-u)\big)_N$ on the set $\{(r,u)\in \R^2; |r-u|\geq d\}$ towards $K^\epsilon$ and the weak convergence of the measure $M^N$ towards $Y|\cdot|$ (cf. Lemma \ref{birk}), the random variable
$$\int_{A\times B}e^{\sum_{n=0}^N\bar{k}_n(r-u)} M^N(dr)M^N(du)$$ almost surely converges towards
$$Y^2\int_{A\times B}e^{K^\epsilon(r-u)} dr\,du. $$
The dominated convergence theorem then yields:
$$\E_Y\big[\min\big(R,M(A)M(B)\big)\big ] \leq \E_Y\Big[\min\Big(R,Y^2\int_{A\times B}e^{K^\epsilon(r-u)} dr\,du\Big)\Big].$$
By letting $R\to \infty$, the monotone convergence theorem yields 
$$\E_Y\big[M(A)M(B)\big ] \leq Y^2\int_{A\times B}e^{K^\epsilon(r-u)} dr\,du.$$
On the other hand, we also have
\begin{equation}\label{ecs2}
\E_Y\big[M(A)M(B)\big ] =\E_Y\big[\widetilde{M}^N(A)\widetilde{M}^N(B)\big]=\E_Y\big[\E\big[\widetilde{M}^N(A)\widetilde{M}^N(B)|\mathcal{G}\big]\big].
\end{equation}
By gathering \eqref{ecs} and \eqref{ecs2} and by using the Fatou's lemma, we deduce
$$\E_Y\big[M(A)M(B)\big ] \geq  Y^2\int_{A\times B}e^{K^\epsilon(r-u)} dr\,du.$$ This completes the proof.\qed

\begin{lemma}\label{cvke}
We have
$$\sup_{|r|\geq d}|K^\epsilon(r)|\to 0 \quad \text{ as }d\to \infty.$$
\end{lemma}

\noindent {\it Proof.} By using \eqref{compk}, we have for $|r|\geq d$:
$$
|K^\epsilon(r)| \leq  \frac{C_\epsilon}{-\ln \epsilon}\int_{\epsilon d}^{+\infty}\theta(u)\ln \frac{u}{\epsilon d}\,du$$
Now, if $\epsilon d\geq 1$, we have:
\begin{align*}
\sup_{|r|\geq d}|K^\epsilon(r)|\leq & \frac{C_\epsilon}{-\ln \epsilon}\int_{\epsilon d}^{+\infty}\theta(u)\ln u\,du
\end{align*}
Hence the result follows from the convergence of the last integral.\qed

\begin{proposition}\label{mixing}
The measure $M$ possesses the following mixing property: given two disjoint sets $A,B$ such that ${\rm dist}(A,B)=d>0$ we have:
\begin{equation}
\big|\E_Y\big[M(A)M(B)]-Y^2|A||B|\big|\leq Y^2\xi(d)|A||B|
\end{equation}
for some function $\xi:\R_+\to\R_+$ such that $\lim_{d\to \infty}\xi(d)=0$.

As a consequence, for any Lebesgue integrable function $\phi$ on $\R^2$ and $d>0$, we have:
\begin{equation}\label{mixfunc}
\Big|\E_Y\big[\int_{|u-r|>d}\phi(u,r)M(dr)M(du)]-Y^2\int_{|u-r|>d}\phi(u,r)\,du\,dr\Big|\leq Y^2\xi(d)\int_{|u-r|>d}|\phi(u,r)|\,du\,dr.
\end{equation}

\end{proposition}

\noindent {\it Proof.} From Lemma \ref{lemmix}, we have
\begin{align*}
\big|\E_Y\big[M(A)M(B)]-Y^2|A||B|\big|&=Y^2\int_{A\times B}(e^{K^\epsilon(r-u)}-1)dr\,du\\
&\leq Y^2\varepsilon(d)|A||B|
\end{align*}
where we have set $\xi(d)=\sup_{|r|\geq d}|e^{K^\epsilon(r)}-1|$. From Lemma \ref{cvke}, we have $\lim_{d\to \infty}\xi(d)=0$. It is then plain to derive \eqref{mixfunc}.
\qed

As a direct consequence, we obtain:
\begin{corollary}\label{coromix}
For any Lebesgue integrable function $\phi$ on $\R^2$ and $d>0$, we have for all $N\in\N\setminus\{0\}$:
\begin{equation*}
\begin{split}
\big|\int_{|u-r|>d}\phi(u,r)\E_Y\big[M^N(dr)M^N(du)]-Y^2\int_{|u-r|>d}|&\phi(u,r)|\,du\,dr\big|\\&\leq Y^2\xi\big(\frac{d}{\epsilon^{N}}\big)\int_{|u-r|>d}|\phi(u,r)|\,du\,dr.
\end{split}
\end{equation*}

\end{corollary}

\subsection{Characterization of the measure $M$}

Having in mind that the measure $M^N$ weakly converges towards $Y|\cdot|$ as $N$ goes to infinity, it is very tantalizing to think that the solution of our problem reduces to taking the limit in \eqref{mn} as $N\to\infty$.
 However, multiplicative chaos badly behaves with respect to weak convergence of measures. So we want to get rid of the measure $M^N$ and have the Lebesgue measure instead in order to deal with a multiplicative chaos in the sense of Kahane. This is the main difficulty of the proof. For that purpose, it is appropriate to take the conditional expectation of $\widetilde{M}^N$ with respect to the $\sigma$-algebra  $\mathcal{F}_N=\sigma(X^0,\dots,X^N,Y)$. Therefore, for any Borelian subset $A$ of $\R$,  we define
$$G_N(A)=\E[\widetilde{M}^N(A)|\mathcal{F}_N]$$ and we claim
\begin{lemma}\label{lem:gn}
The following relation holds for each $N\geq 0$:
\begin{equation}\label{eq:birk}
G_N(A)=Y\int_A\exp\Big(\sum_{n=0}^NX^n_r-\frac{1}{2}\E[(X^n_r)^2]\Big)\,dr.
\end{equation}
Furthermore, for each bounded Borelian set $A$, the sequence $(G_N(A))_N$ is a  positive  martingale bounded in $L^{1+\delta}$.
\end{lemma}

\noindent {\it Proof.} If $A$ has infinite Lebesgue measure, both sides of \eqref{eq:birk} are infinite. So we focus on the case when $A$ has finite Lebesgue measure. First observe that for each $s<t$ and $A\in \mathcal{F}_N$, we have from Lemma \ref{birk}
$$\E[\int_\R\ind_{[s,t]}(r)\ind_A M^N(dr)|\mathcal{F}_N] =\ind_A\E_Y[M^N([s,t])]=\ind_AY(t-s).$$
By using density arguments and Fatou's lemma, we establish that, for each positive $\mathcal{F}_N\otimes \mathcal{B}(\R)$-measurable function $\varphi\in L^1(\Omega\times \R;\p\otimes dt)$, we have
$$\E\big[\int_\R\varphi(\omega,r)M^N(dr)\big|\mathcal{F}_N]=\int_\R\varphi(\omega,r)Y\,dr.$$
So \eqref{eq:birk} is proved.

 Finally, for each bounded set $A$ we have $\E[M(A)^{1+\delta}]<+\infty$ for some $\delta>0$. The Jensen inequality then yields
 \begin{align*}
 \E[(G_N(A))^{1+\delta}]&= \E[(\E[\widetilde{M}^N(A)|\mathcal{F}_N])^{1+\delta}]\leq \E[(\widetilde{M}^N(A))^{1+\delta}] =\E[M(A)^{1+\delta}]<+\infty.
 \end{align*}
The martingale $(G_N(A))_N$ is thus bounded in $L^{1+\delta}$.\qed

\vspace{2mm}
Being bounded in $L^{1+\delta}$, the martingale converges almost surely and in $L^{1+\delta}$ towards a random variable $Q(A)$, which can be formally thought of as 
$$Q(A)= Y \int_A\exp\left( X_r-\frac{1}{2} \E[X_r^2] \right) \,dr$$ where $\left(X_r\right)_{r \in \R}$ is a "Gaussian process" with covariance kernel $K^\epsilon(r)$, that is a Gaussian multiplicative chaos. The remaining part of our argument can be roughly  summed up as follows. First, we obtain estimates on the kernel $K^\epsilon$ derived from the fact that the Gaussian multiplicative chaos $Q$ admits a moment of order $1+\delta$. Second, we use these estimates to prove that $Q $ has  the same law as $M$. Finally, since $Q$ has the same law as $M$, which does not depend on $\epsilon$, the kernel $K^\epsilon$ should not depend on $\epsilon$ either. This is a strong constraint on $K^\epsilon$, from which we derive the specific structure of $K^\epsilon$ given by \eqref{struct}. 

So we claim
\begin{proposition}\label{estimate}
For each $0<\gamma<\delta$, we can find $\rho>0$ such that:
\begin{equation}
\sup_{n}n^{1+\rho}\E[M([0,\frac{1}{n}])^{1+\gamma}]<+\infty.
\end{equation}
\end{proposition}

\noindent {\it Proof.} The proof relies on the following bound (see the proof below):

\begin{lemma}\label{bound}
The existence of a moment of order $1+\delta$ for the measure $M$ implies the following bound:
$$ k_\epsilon(0)\leq \frac{2}{1+\delta}\ln \frac{1}{\epsilon}.$$
\end{lemma}

Since we have for all $r\in\R$: $k_\epsilon(r)\leq k_\epsilon(0)$, the covariance kernel of the process $\omega_\epsilon$ is dominated by that of the constant process $\omega_\epsilon(0)$. Hence, by using \eqref{star} and Lemma \ref{cvx}, it is plain to see that, for each $\gamma>0$: 
\begin{align*}
\E[M([0,\frac{1}{n}])^{1+\gamma}]&=\E\left[\Big(\int_0^{1/n}e^{\omega_{1/n}(r)}M^{1/n}(dr)\Big)^{1+\gamma}\right]\\
&\leq \E\left[\Big(\int_0^{1/n}e^{\omega_{1/n}(0)}M^{1/n}(dr)\Big)^{1+\gamma}\right]\\
&\leq \E\left[e^{(1+\gamma)\omega_{1/n}(0)}\right]\E\left[\Big(M^{1/n}([0,\frac{1}{n}])\Big)^{1+\gamma}\right]\\
&=e^{\frac{(1+\gamma)^2}{2}k_{1/n}(0)-\frac{1+\gamma}{2}k_{1/n}(0)}\E\Big[\Big(M([0,1])\Big)^{1+\gamma}\Big]\frac{1}{n^{1+\gamma}}.
\end{align*}
Since $k_{1/n}(0)\leq \frac{2}{1+\delta}\ln n$, we deduce
\begin{align*}
\E[M([0,\frac{1}{n}])^{1+\gamma}]\leq & e^{\big(\frac{\gamma^2+\gamma}{1+\delta}-\gamma-1\big)\ln n}\E\Big[\Big(M([0,1])\Big)^{1+\gamma}\Big]\\
=& \frac{1}{n^{1+\rho}}\E\Big[\Big(M([0,1])\Big)^{1+\gamma}\Big]
\end{align*}
where we have set $$\rho\stackrel{def}{=}-\frac{\gamma^2+\gamma}{1+\delta}+\gamma.$$ Clearly, we have $\rho>0$ provided that $0<\gamma<\delta $. 
The proof of Proposition \ref{estimate} is complete.\qed

\vspace{3mm}
\noindent {\it Proof of Lemma \ref{bound}.} Let $n \in \N$.
\begin{align}
\E\left[M[0;t]^{1+\delta}\right] &= \E\left[\left( M[0;\frac{t}{n}] +  M[\frac{t}{n};\frac{2t}{n}] + \dots + M[\frac{(n-1)t}{n}; t] \right)^{1+\delta}\right]  \\
                                                  &\geq  \E\left[\left( M[0;\frac{t}{n}] \right)^{1+\delta} + \left( M[\frac{t}{n};\frac{2t}{n}] \right)^{1+\delta} + \dots 
+ \left( M[\frac{(n-1)t}{n}; t] \right)^{1+\delta}  \right] \\ 
&= n  \E\left[\left( M[0;\frac{t}{n}] \right)^{1+\delta} \right]
\label{stationarity}
\end{align}
We used the stationarity of the measure $M$ in the second line.
Now write, for $h>0$:
\begin{equation}
g(h) = \sup_{r\leq h} \mid k_{1/n}(0) - k_{1/n}(r) \mid
\end{equation}
We have, for every $r \in (0,t/n]$ and $n$ large enough: 
\begin{equation*}
\mid k^{1/n}(0) - g(t/n) \mid \leq k^{1/n}(r).
\end{equation*}
So, using classical Gaussian inequality (see Lemma \ref{cvx}): 
\begin{align}
\E\left[M[0;\frac{t}{n}]^{1+\delta}\right] &= \E\left[\left( \int_0^{t/n} e^{\omega_{1/n}(r)} M^{1/n}(dr) \right)^{1+\delta}\right]  \nonumber\\
  &\geq  \E\left[\left(   \int_0^{t/n} e^{\sqrt{\mid k_{1/n}(0) - g(t/n) \mid} Z_n - \frac{1}{2} \mid k_{1/n}(0) - g(t/n) \mid} M^{1/n}(dr) \right)^{1+\delta}\right]\nonumber \\
  &= \E\left[\left( e^{\sqrt{\mid k_{1/n}(0) - g(t/n) \mid} Z_n -\frac{1}{2} \mid k_{1/n}(0) - g(t/n) \mid} \right)^{1+\delta}\right] \E\left[\left(M^{1/n}[0;\frac{t}{n}]\right)^{1+\delta}\right] \nonumber\\ 
  &= e^{-\frac{1+\delta}{2}  \mid k_{1/n}(0) - g(t/n) \mid} e^{\frac{(1+\delta)^2}{2} \mid k_{1/n}(0) - g(t/n) \mid} \frac{1}{n^{1+\delta}} \E\left[\left(M[0;t]\right)^{1+\delta}\right] \label{JensenStareq}
\end{align}
We used Lemma $13$ in the second line.
Using equations \eqref{stationarity} and \eqref{JensenStareq}, one gets 
\begin{equation} \label{eqwithg}
 e^{-\frac{1+\delta}{2}  \mid k_{1/n}(0) - g(t/n) \mid} e^{\frac{(1+\delta)^2}{2} \mid k_{1/n}(0) - g(t/n) \mid} \frac{1}{n^{\delta}} \leq 1
\end{equation}
As $h$ goes to $0$, $g(h)$ goes to $0$ (the function $k_{1/n}$ is continuous). Letting $t$ goes to $0$ in (\ref{eqwithg}), one gets 
\begin{equation*}
k_{1/n}(0) \leq \frac{2}{1+\delta} \ln n.
\end{equation*}
and the lemma is proved. \qed

\vspace{2mm}
We are now in position to tackle the main step of the proof:
\begin{proposition}\label{law}
The random measures $(Q(A))_{A\in \mathcal{B}(\R)}$ and $(M(A))_{A\in \mathcal{B}(\R)}$ have the same law.
\end{proposition}

\noindent {\it Proof.} Let $F$ be some function defined on $\R_+$ such that:
\begin{itemize}
\item $F$ is convex,
\item $F(x)\leq Cx^{1+\gamma}$ for some constants $C>0$ and $0<\gamma<\delta$,
\item $F\circ \sqrt{\phantom{5}}$ is concave, nondecreasing and sub-additive.
\end{itemize}
Let $f$ be a lower semi-continuous positive function on $\R$ with compact support. We have by Jensen's inequality: 
\begin{align*}
\E\Big[F\big(\int_\R f(x)\,M(dx)\big)\Big]&=\E\Big[F\big(\int_\R f(x)\,\widetilde{M}^N(dx)\big)\Big]\\
&=\E\Big[\E\Big[F\big(\int_\R f(x)\,\widetilde{M}^N(dx)\big)|\mathcal{F}_N\Big]\Big]\\
&\geq\E\Big[F\big(\int_\R f(x)\,G_N(dx)\big)\Big].
\end{align*}
We let $N$ go to $+\infty$. By using the weak convergence of  $G_N(dr)$ towards $Q(dr)$, we obtain:
\begin{equation}
\E\big[F\big(\int_\R f(r)\,M(dr)\big)\big]\geq \E\big[F\big(\int_\R f(r)\,Q(dr)\big)\big].
\end{equation}

Now we want to establish the converse inequality. We set $\widetilde{F}=F\circ \sqrt{\phantom{5}}$.  For any $\tau>0$, we have by using  the sub-additivity of $\widetilde{F}$: 
\begin{align}
\E\big[F\big(\int_\R f(r)\,M(dr)\big)\big]=&\E\Big[\widetilde{F}\Big(\big(\int_\R f(r)\,\widetilde{M}^N(dr)\big)^2\Big)\Big]\nonumber\\
=&\E\Big[\widetilde{F}\Big(\int_\R\int_\R f(r)f(u)\,\widetilde{M}^N(dr)\widetilde{M}^N(du)\Big)\Big]\nonumber\\
\leq &\E\Big[\widetilde{F}\Big(\int_{|r-u|\leq \tau} f(r)f(u)\,\widetilde{M}^N(dr)\widetilde{M}^N(du)\Big)\Big]\nonumber\\
&+\E\Big[\widetilde{F}\Big(\int_{|r-u|> \tau} f(r)f(u)\,\widetilde{M}^N(dr)\widetilde{M}^N(du)\Big)\Big]
\nonumber.
\end{align}
Then, by conditioning with respect to $\mathcal{F}_N$ and by using the Jensen inequality in the second term of the latter inequality, we deduce:
\begin{align}
\E\big[F&\big(\int_\R f(r)\,M(dr)\big)\big]\\
\leq &\E\Big[\widetilde{F}\Big(\int_{|r-u|\leq \tau} f(r)f(u)\,\widetilde{M}^N(dr)\widetilde{M}^N(du)\Big)\Big]\nonumber\\
&+\E\Big[\widetilde{F}\Big( \int_{|r-u|> \tau} f(r)f(u)\exp\big(\sum_{k=0}^NX^n_r+X^n_u-k_n(0)\big)\,\E_Y[M^N(dr)M^N(du)]\Big) \Big]\nonumber\\
 \stackrel{def}{=} & C(1,\tau,N)+C(2,\tau,N)\label{c12N}.
\end{align}

We claim: 
\begin{lemma}\label{c2N}
For each fixed $\tau>0$, $C(2,\tau,N)$ converges as $N\to \infty$ towards 
\begin{equation*}
\E\Big[\widetilde{F}\Big( \int_{|r-u|> \tau} f(r)f(u) \,Q(dr)Q(du)\Big) \Big].
\end{equation*}
Furthermore, this latter quantity converges, as $\tau\to 0$, towards
\begin{equation*}
\E\Big[F\Big( \int f(r) \,Q(dr) \Big) \Big].
\end{equation*}
Finally, the quantity $C(1,\tau,N)$ converges to $0$ as $\tau\to 0$ uniformly with respect to $N\in\N^*$.
\end{lemma}

Let us admit for a while the above lemma to finish the proof of Proposition \ref{law}. By gathering \eqref{c12N} and  Lemma \ref{c2N}, we deduce
\begin{align*}
\E\big[F\big(\int_\R f(r)\,M(dr)\big)\big]\leq & \liminf_{\tau\to 0}\E\Big[\widetilde{F}\Big( \int_{|r-u|> \tau} f(r)f(u) \,Q(dr)Q(du)\Big) \Big]\\
=&\E\Big[F\Big( \int f(r) \,Q(dr) \Big) \Big] .
\end{align*}
Hence we have proved
\begin{equation}
\E\big[F\big(\int_\R f(r)\,M(dr)\big)\big]=\E\big[F\big(\int_\R f(r)\,Q(dr)\big)\big]. 
\end{equation}

The basic choice for $F$ is the function $x\mapsto x^{1+\gamma}$ with $0<\gamma<\delta$. Thus we have proved that the mappings
$$\E\big[\exp\big(z\ln\int_\R f(r)\,M(dr)\big)\big]\quad \text{ and }\quad\E\big[\exp\big(z\ln \int_\R f(r)\,Q(dr)\big)\big] $$ 
coincide  for $z \in ]1,1+\delta[$.
By analyticity arguments, we deduce that $\int_\R f(x)\,M(dx)$  and $\int_\R f(x)\,Q(dx)$ have the same law. This is enough to prove that the random measures $M$ and $Q$ have the same law. Indeed, if we consider two families $(\lambda_i)_{1\leq i\leq n}$ of positive real numbers and $(A_i)_{1\leq i\leq n}$ of bounded open subsets  of $\R$, we  define the lower semi-continuous function
$$f(x)=\sum_{i=1}^n\lambda_i\ind_{A_i}(x)$$ and we obtain
$$\sum_{i=1}^n \lambda_iM(A_i)\stackrel{law}{=}\sum_{i=1}^n \lambda_iQ(A_i).$$ 
 It turns out that the law of a random vector $(Y_1,\dots,Y_n)$ made up of positive random variables is characterized by the  combinations 
$$\sum_{i=1}^n\lambda_iY_i$$ where $(\lambda_i)_{1\leq i\leq n}$ is a family of positive real numbers. The proof of Proposition \ref{law} is complete.\qed

\vspace{2mm}
{\it Proof of Lemma \ref{c2N}.} Let us first investigate the quantity $C(1,\tau,N)$. Assume the function $f$ has its support included in the ball $B(0,R)$ for some $R>0$. We can cover the set $$\{(x,y)\in \R^2;|x-y|\leq \tau\text{ and }\max(|x|,|y|)\leq R\}$$ by the squares 
$$A^n_j=[t_j^n,t_{j+2}^n]\times [t_j^n,t_{j+2}^n]\quad \text{where }t_j^n=-R+2\tau j, \,\,\text{ for }j=0,\dots,E(\frac{R}{\tau}).$$  
We set $S=\sup_{ \R} f$. Because  $\widetilde{F}$ is sub-additive and increasing, we have: 
\begin{align*}
C(1,\tau,N)&\leq  \E\Big[\widetilde{F}\Big(\sum_{0\leq j\leq E(\frac{R}{\tau})}\int_{A^n_j} f(r)f(u)\,\widetilde{M}^N(dr)\widetilde{M}^N(du)\Big)\Big]\\
&\leq \sum_{0\leq j\leq E(\frac{R}{\tau})} \E\Big[\widetilde{F}\Big(\int_{A^n_j} f(r)f(u)\,\widetilde{M}^N(dr)\widetilde{M}^N(du)\Big)\Big]\\
&\leq \sum_{0\leq j\leq E(\frac{R}{\tau})} \E\Big[\widetilde{F}\Big(S^2\int_{A^n_j} \,\widetilde{M}^N(dr)\widetilde{M}^N(du)\Big)\Big]\\
&= \sum_{0\leq j\leq E(\frac{R}{\tau})} \E\Big[\widetilde{F}\Big(S^2(\widetilde{M}^N([t_j^n,t_{j+2}^n]))^2\Big)\Big]\\
&= \sum_{0\leq j\leq E(\frac{R}{\tau})} \E\Big[F\Big(SM([t_j^n,t_{j+2}^n])\Big)\Big].
\end{align*}
By stationarity, we deduce
\begin{align*}
C(1,\tau,N)&\leq \frac{2R}{\tau}\E\Big[F\Big(SM([0,2\tau])\Big)\Big]\\
&\leq  \frac{2R}{\tau}S^{1+\gamma}\E\big[M([0,2\tau])^{1+\gamma}\big].
\end{align*}
It results from Proposition  \ref{estimate} that the last quantity
converges towards $0$ as $\tau$ goes to $0$ uniformly with respect to $N$.

Now we investigate the quantity $C(2,\tau,N)$. Since $\tilde{F}$ is sub-additive and increasing, we have $|\tilde{F}(a)-\tilde{F}(b)|\leq \tilde{F}(|b-a|)$ for all positive real numbers $a,b$. This together with Corollary \ref{coromix} yields
\begin{align*}
\Big|C(2,\tau,N)-&\E\Big[\widetilde{F}\Big( \int_{|r-u|> \tau} f(r)f(u)\exp\big(\sum_{k=0}^NX^n_r+X^n_u-k_n(0)\big)\,Y^2dr\,du\Big) \Big]\Big|\\
\leq &\E\Big[\widetilde{F}\Big(Y^2\xi \big(\frac{\tau}{\epsilon^N}\big)\int_{|r-u|> \tau} f(r)f(u)\exp\big(\sum_{k=0}^NX^n_r+X^n_u-k_n(0)\big)\,dr\,du \Big) \Big]\\
\leq &\E\Big[\widetilde{F}\Big(\xi \big(\frac{\tau}{\epsilon^N}\big)S^2G_N([-R,R])^2 \Big)  \Big]\\
\leq &\E\Big[F\Big(S\xi \big(\frac{\tau}{\epsilon^N}\big)^{1/2}G_N([-R,R])  \Big) \Big]\\
\leq &\xi \big(\frac{\tau}{\epsilon^N}\big)^{\frac{1+\gamma}{2}}S^{1+\gamma} \E\Big[G_N([-R,R]) ^{1+\gamma}\Big].
\end{align*}
Obviously, the last quantity converges to $0$ as $N$ goes to $\infty$. Furthermore, the quantity $$\widetilde{F}\Big( \int_{|r-u|> \tau} f(r)f(u)\exp\big(\sum_{k=0}^NX^n_r+X^n_u-k_n(0)\big)\,Y^2dr\,du\Big) $$ almost surely converges towards $$\widetilde{F}\Big( \int_{|r-u|> \tau} f(r)f(u)\,Q(dr)\,Q(du)\Big) $$ and is uniformly integrable because $F(x)\leq Cx^{1+\gamma}$ and $Q$ is a multiplicative chaos admitting a moment of order $1+\delta$ with $\delta>\gamma$. The Lebesgue convergence theorem then yields:
\begin{align*}
\E\Big[\widetilde{F}\Big(& \int_{|r-u|> \tau} f(r)f(u)\exp\big(\sum_{k=0}^NX^n_r+X^n_u-k_n(0)\big)\,Y^2dr\,du\Big) \Big]\\
&\to \E\Big[\widetilde{F}\Big( \int_{|r-u|> \tau} f(r)f(u)\,Q(dr)\,Q(du)\Big)\Big]\quad \text{ as } N\to \infty. 
\end{align*}
Gathering the above relations yields
$$C(2,\tau,N)\to \E\Big[\widetilde{F}\Big( \int_{|r-u|> \tau} f(r)f(u)\,Q(dr)\,Q(du)\Big)\Big]\quad \text{ as } N\to \infty. $$

Similar arguments as those used above allow to establish that 
\begin{align*}
\liminf_{\tau\to 0}\E\Big[\widetilde{F}\Big( \int_{|r-u|> \tau} f(r)f(u) \,Q(dr)Q(du)\Big) \Big]=&
\E\Big[\widetilde{F}\Big( \int_{\R^2} f(r)f(u) \,Q(dr)Q(du)\Big) \Big]\\
=&\E\Big[F\Big( \int_{\R} f(r) \,Q(dr)\Big) \Big].
\end{align*}
 Indeed, by proceeding as for $C(1,\tau,N) $, we can prove that the "diagonal contribution" goes to $0$ as $\tau\to 0$. Details are left to the reader.
The proof of the Lemma is complete. \qed

 \vspace{2mm}
The final step of our argument is now to prove that the kernel $K^\epsilon$ defined by \eqref{ke} does not depend on $\epsilon$. Expressing the kernel $K^\epsilon$ as a function of the marginals of the measure $M$ is enough for that purpose. So we remind the reader of Lemma \ref{lemmix}, which states
$$\E_Y[M(A)M(B)]=Y^2\int_{A\times B}e^{K^\epsilon(r-u)}\,drdu.$$
We deduce that, for any $s\not =0$ and on the set $\{Y>0\}$,
\begin{equation}\label{limh}
K^\epsilon(s)=\lim_{h\to 0}\ln\Big(\frac{1}{h^2}\E_Y[M([0,h])M([s,s+h])]\Big)-2\ln Y.
\end{equation}

As a straightforward consequence, the kernel $K^\epsilon$ defined by \eqref{ke} does not depend on $\epsilon$ since the left-hand side in \eqref{limh} does not either. So we can define the quantity
$$\forall r\not =0,\quad K(r)=K^\epsilon(r)$$ for some $\epsilon\in (0,1)$ and this relation is also valid for any $\epsilon\in (0,1)$. It is also plain to see that for each $\epsilon\in (0,1)$ we have:
\begin{equation}\label{relK}
\forall r\not =0,\quad K(r)=k_\epsilon(r)+K(\frac{r}{\epsilon})
\end{equation}
since $K^\epsilon$ satisfies such a relation. Such a specific  functional equation  implies a precise structure for the function $K$:

\begin{proposition}\label{formK}
For $r>0$, we have 
\begin{equation}
K(r) = \int_r^{+\infty} \frac{k(u)}{u} du
\end{equation}
where $k(u)$ is a positive-definite continuous function $\R_+ \to \R$.
\end{proposition}

\noindent{\it{Proof.}} Because  $K$ is Lipschitzian on the compact subsets of $\R\setminus\{0\}$, there exists a locally bounded
measurable function $f$ on $(0;+\infty)$ such that for all $r,s> 0$, 
\begin{equation*}
K(s) - K(r)  = \int_r^s f(t) dt. 
\end{equation*}
Define, for $r\in \R$, 
\begin{equation*}
\phi(r) = K(e^r)
\end{equation*}
It is straightforward to derive from \eqref{relK} that, for all $r\in \R, \alpha \geq 0 $, 
\begin{equation} \label{derive}
\phi(r+\alpha) - \phi(r) = - k_{e^{-\alpha}}(e^r)
\end{equation}
Note that $ k_{1}(e^r) = 0$. From equation (\ref{derive}), one obtains :
\begin{equation} \label{derive2}
 \frac1\alpha \int_r^{r+\alpha} e^u f(e^u) du = - \frac{k_{e^{-\alpha}}(e^r)}{\alpha}
\end{equation}

For almost every $r$, the left-hand side of equation (\ref{derive2}) tends to $e^r f(e^r)$ when $\alpha$ goes to $0$. 
Thus, the right-hand side of $(\ref{derive2})$ converges also for almost every $r$ to $e^r f(e^r)$ when  $\alpha$ goes to $0$. 
 
We define the function $g$ by the following limit for almost every $r$:
\begin{equation}
g(r)=\underset{\alpha \to 0}{\lim} -\frac{1}{\alpha}\int_{r}^{r+\alpha}e^{u}f(e^u)du= \underset{\alpha \to 0}{\lim} \frac{k_{e^{-\alpha}}(e^r)}{\alpha}
\end{equation}

As defined, the function $g$ is measurable with 
respect to the Borelian $\sigma$-field of $\R$.
For almost every $x \in (0,+\infty)$, define $$h(x) = g(\ln(x)),$$
and $h(0)$ by $h(0)=\frac{k_{e^{-\alpha}}(0)}{\alpha}$ for some $\alpha>0$. Note that the definition of $h(0)$ does not depend on $\alpha$ because:
\begin{lemma}\label{form}
We have the following asymptotic behaviour of $K$ around $0$:
$$K(r)\simeq \frac{k_\epsilon(0)}{\ln \epsilon}\ln r\quad \text{ as }r\to 0.$$
\end{lemma}

Thus $h$ is well defined at $0$ and we can now prove that it is positive definite:
\begin{lemma}
The function $h(|.|)$ is positive definite (as a tempered distribution in the sense of Schwartz, see \cite{cf:Gel} or \cite{cf:Schw}). One can also find a symmetric positive measure $\mu$ on $\R$ (with $\mu(\R)<\infty$) such that for almost every $x \in \R$:   
\begin{equation*}
h( |x| ) = \int_\R e^{ix\xi} \mu(d\xi)
\end{equation*} 
\end{lemma}

\noindent{\it Proof.}
For almost every $x \in \R$, $h(|x|)=\underset{\alpha \to 0}{\lim} \frac{k_{e^{-\alpha}}(|x|)}{\alpha}$ and $\frac{k_{e^{-\alpha}}(|x|)}{\alpha} \leq h(0)$ uniformly in $\alpha$. Thus, if $\varphi$ is a smooth function with compact support, we get using the dominated convergence theorem:
\begin{equation*}
 \int_{\R}  \int_{\R} h( |y-x| ) \varphi(x)\overline{\varphi}(y)dxdx  = \underset {\alpha \to 0} {\lim} \int_{\R}  \int_{\R} \frac{k_{e^{-\alpha}}(|y-x|)}{\alpha} \varphi(x)\overline{\varphi}(y)dxdx \geq 0.
\end{equation*}
We conclude that $h(|.|)$ is positive definite. By the Bochner-Schwartz theorem, the Fourier transform of $h(|.|)$ is a symmetric positive measure $\mu(d\xi)$ such that there exists $p \geq 0$ with:
\begin{equation*}   
\int_{\R}\frac{\mu(d\xi)}{(1+|\xi|)^p} < \infty.
\end{equation*}
In order to conclude, it is sufficient to prove that $\mu(\R)< \infty$. We note $\theta(x)=\frac{e^{-x^2/2}}{\sqrt{2 \pi}}$ and $\theta^{\epsilon}=\frac{1}{\epsilon} \theta(./ \epsilon)$ for $\epsilon >0$. By the inverse Fourier theorem, we get:
\begin{equation*}
(\theta^{\epsilon} \ast h)(0)=\int_{\R} e^{-\epsilon^2 \xi^2 /2} \mu(d\xi).
\end{equation*}
Thus the right hand side of the above equality is bounded by $h(0)$ and we conclude by letting $\epsilon$ go to $0$. \qed

Integrating with respect to the Lebesgue measure the relation $g(t) = - e^t f(e^t)$ which is true for almost every $t\in \R$, one gets
\begin{equation*}
K(s) - K(r) = - \int_r^{s} \frac{h(u)}{u} du.
\end{equation*} 
Because $K(s)\to 0$ as $s\to +\infty$, the function $u\mapsto \frac{h(u)}{u}$ is integrable at the vicinity of $+\infty$ in the generalized sense.
We deduce:
\begin{equation*}
K(r) = \int_r^{+\infty} \frac{h(u)}{u} du.
\end{equation*} 

By the previous lemma, there exists a finite symmetric positive measure $\mu$ on $\R$ such that, for almost every $x \in \R$, 
\begin{equation*}
h(x) = \int_\R e^{ix\xi} \mu(d\xi)
\end{equation*} 
For simplicity, define for all $x \in \R$, $k(x) = \int_\R e^{ix\xi} \mu(d\xi)$. The function $k$ is continuous on $\R$. We get finally, 
\begin{equation}
K(r) = \int_r^{+\infty} \frac{k(u)}{u} du.
\end{equation} 
The proof of Proposition \ref{formK} is complete.\qed

\vspace{2mm}
\noindent {\it Proof of Proposition \ref{propcut}.} This is just a direct consequence of Theorem \ref{main} and equation \eqref{limh}.\qed

\appendix
\section{Proofs of some auxiliary lemmas}\label{app:proof}
 
 \begin{lemma}\label{lemY}
Let $F : \R^n \mapsto \R$ be a measurable function. Then, for all bounded Borelian sets $A_1,\dots,A_n \subset \R$, the following relation holds almost surely:
\begin{equation*}
\E_Y\left[F(M(A),\dots,M(A_n))\right] = \E_Y\left[F(\tilde{M}^N(A),\cdots,\tilde{M}^N(A_n))\right]
\end{equation*}
\end{lemma}

\noindent{\it{Proof.}} By using the Jensen inequality, we have 
\begin{align*}
\E\Big[\Big|&\frac1T \tilde{M}^N\left[0;T\right] -\frac1T M^N\left[0;T\right]\Big|\Big]\\
 =&\E\Big[\Big(\Big|\frac1T \tilde{M}^N\left[0;T\right] -\frac1T M^N\left[0;T\right]\Big|^2\Big)^{1/2}\Big]\\
\leq & \E\Big[\Big(\E\Big[\Big|\frac1T \tilde{M}^N\left[0;T\right] -\frac1T M^N\left[0;T\right]\Big|^2|M\Big]\Big)^{1/2}\Big]\\
= & \E\Big[\Big(\frac{1}{T^2}\int_0^T\int_0^T\E\Big[\big(e^{\sum_{n=0}^NX^n_r-\frac{1}{2}\E[(X^n_r)^2]}-1\big)\big(e^{\sum_{n=0}^NX^n_u-\frac{1}{2}\E[(X^n_u)^2]}-1\big)\Big]M^N(dr)M^N(du)\Big)^{1/2}\Big]\\
=& \E\Big[\Big(\frac{1}{T^2}\int_0^T\int_0^T\big(e^{\sum_{n=0}^N\bar{k}_n(r-u)}-1\big) M^N(dr)M^N(du)\Big)^{1/2}\Big]\\
\end{align*}
The integrand in the above expectation converges almost surely towards $0$ because, for each $0\leq n\leq N$,  $\bar{k}_n$ is bounded and converges to $0$ in the vicinity of $\infty$. Furthermore, it is uniformly integrable because
$$\sup_T\E\Big[\Big(\frac1T M^N(\left[0;T\right])\Big)^{1+\delta}\Big]<+\infty.$$
We deduce that 
$$\E\Big[\Big|\frac1T \tilde{M}^N\left[0;T\right] -\frac1T M^N\left[0;T\right]\Big|\Big]\to 0 \quad \text{ as }T\to +\infty.$$
As a consequence, $\frac1T \tilde{M}^N\left[0;T\right]$ converges  almost surely along a subsequence towards $Y$.

One has, for any function $h$ bounded and continuous,
\begin{equation*}
\E\left[F(M(A_1),\dots,M(A_n))h\left(\frac1T M\left[0;T\right] \right)\right] = \E\left[F(\tilde{M}(A_1),\dots,\tilde{M}(A_n))h\left(\frac1T \tilde{M}^N\left[0;T\right]\right)\right]
\end{equation*} 
Sending $T$ to $+\infty$ along the subsequence, we get by the bounded convergence theorem
\begin{equation*}
\E\left[F(M(A_1),\dots,M(A_n)) h\left(Y\right)\right] = \E\left[F(\tilde{M}(A_1),\dots,\tilde{M}(A_n)) h\left(Y\right)\right]
\end{equation*} 
and the lemma is proved. \qed

\begin{lemma}\label{cvx}
Let $F:\R_+\to \R$ be some convex  function such that 
$$\forall x\in\R_+,\quad |F(x)|\leq M(1+|x|^\beta),$$ for some positive constants $M,\beta$, and $\sigma$ be a Radon measure on the Borelian subsets of $\R$. Given $a<b$, let $(X_r)_{a\leq r\leq b},(Y_r)_{a\leq r\leq b}$ be two continuous centered Gaussian processes with continuous covariance kernels $k_X$ and $k_Y$ such that
$$\forall u,v\in [a,b],\quad k_X(u,v)\leq k_Y(u,v). $$
Then
$$\E\Big[F\Big(\int_a^be^{X_r-\frac{1}{2}\E[X_r^2] }\,\sigma(dr)\Big)\Big]\leq \E\Big[F\Big(\int_a^be^{Y_r-\frac{1}{2}\E[Y_r^2] }\,\sigma(dr)\Big)\Big].$$
\end{lemma}

\noindent {\it Proof.} For each $N\in \N$, we define the smooth subdivision $t^N_p=a+p\frac{b-a}{N}$, $p=0,\dots,N$, of the interval $[a,b]$. We also introduce the random variables 
$$S^X_N=\sum_{p=0}^{N-1}e^{X_{t^N_p}-\frac{1}{2}\E[X_{t^N_p}^2]}\sigma([t^N_p,t^N_{p+1})) \quad \text{and}\quad S^Y_N=\sum_{p=0}^{N-1}e^{Y_{t^N_p}-\frac{1}{2}\E[Y_{t^N_p}^2]}\sigma([t^N_p,t^N_{p+1})) .$$
By classical Gaussian inequalities (see \cite[corollary 6.2]{cf:RoVa} for instance), we have
 $$\forall N\geq 1,\quad \E\Big[F\Big(S^X_N\Big)\Big]\leq \E\Big[F\Big(S^Y_N\Big)\Big].$$
So it just remains to pass to the limit as $N\to \infty$ by using the dominated convergence theorem. By continuity of the processes $X,Y$ the random variables $S^X_N,S^Y_N$  converge almost surely respectively towards $\int_a^be^{X_r-\frac{1}{2}\E[X_r^2] }\,\sigma(dr),\int_a^be^{Y_r-\frac{1}{2}\E[Y_r^2] }\,\sigma(dr)$. Clearly, we have:
$$|F(S^X_N)|\leq M\big(1+|S^X_N|^\beta\big),$$ so that we just have to prove that $|S^X_N|^\beta$ is uniformly integrable (the same argument holds for $|S^Y_N|^\beta$). It is enough to establish that for each $d\in\N$,
$$\sup_N\E\big[(S^X_N)^d\big]<+\infty.$$ We have
\begin{align*}
\E\big[(S^X_N)^d\big]=&\E\Big[\Big(\sum_{p=0}^{N-1}e^{X_{t^N_p}-\frac{1}{2}\E[X_{t^N_p}^2]}\sigma([t^N_p,t^N_{p+1}))\Big)^d\Big]\\
=&\sum_{p_1,\dots,p_d=0}^{N-1}\E\Big[e^{X_{t^N_{p_1}}+\cdots+X_{t^N_{p_d}}}\Big]e^{-\frac{1}{2}(\E[X_{t^N_{p_1}}^2]+\cdots+\E[X_{t^N_{p_d}}^2])}\sigma([t^N_{p_1},t^N_{p_1+1}))\times \cdots\times \sigma([t^N_{p_d},t^N_{p_d+1})) \\
=&\sum_{p_1,\dots,p_d=0}^{N-1}e^{\frac{1}{2}\sum_{i,j=1}^d k_X(t^N_{p_i},t^N_{p_j})}e^{-\frac{1}{2}(\E[X_{t^N_{p_1}}^2]+\cdots+\E[X_{t^N_{p_d}}^2])}\sigma([t^N_{p_1},t^N_{p_1+1}))\times \cdots\times \sigma([t^N_{p_d},t^N_{p_d+1})) \\
\to & \int_a^b\dots\int_a^be^{\frac{1}{2}\sum_{i \not = j}^dk_X(u_i,u_j)}\sigma(du_1)\cdots\sigma(du_d)
\end{align*}
as $N\to \infty$. This completes the proof.\qed

%
%
%

\vspace{3mm}
\noindent {\it Proof of Lemma \ref{form}.} We choose any $\epsilon<1$ and consider $|r|\leq 1$.  Since $k_\epsilon$ is continuous at $0$, we can find, for $\alpha>0$, some $\eta>0$  such that $k_\epsilon(0)-\alpha\leq k_\epsilon(u)\leq k_\epsilon(0) $ for $|u|\leq \eta$. Then we decompose $K$ as
\begin{align*}
K^\epsilon(r)&=\sum_{n=0}^{+\infty}k_\epsilon(\frac{r}{\epsilon^n})\\
&=\sum_{n=0}^{\frac{\ln \frac{r}{\eta}}{\ln \epsilon}-1}  k_\epsilon\big(\frac{r}{\epsilon^n}\big)+\sum_{n=\frac{\ln \frac{r}{\eta}}{\ln \epsilon}}^{+\infty}  k_\epsilon\big(\frac{r}{\epsilon^n}\big)\\
&\stackrel{def}{=}\sum_{n=0}^{\frac{\ln \frac{r}{\eta}}{\ln \epsilon}-1}  k_\epsilon\big(\frac{r}{\epsilon^n}\big)+g_\epsilon(r)
\end{align*}
Let us prove that $g_\epsilon$ is bounded over a neighborhood of $0$. By using \eqref{modulus} and following the computations of \eqref{compk}, we have for $p\in \N$:
$$\sum_{n=p}^{+\infty} | k_\epsilon\big(\frac{r}{\epsilon^n}\big)|\leq \frac{2 C_\epsilon}{-\ln \epsilon}\int_{\frac{r}{\epsilon^{p-1}}}^{+\infty}\theta(u)\ln u\,du.$$
We deduce by taking $p=\frac{\ln \frac{r}{\eta}}{\ln \epsilon}$:
$$|g_\epsilon(r)|\leq  \frac{2 C_\epsilon}{-\ln \epsilon}\int_{\eta\epsilon}^{+\infty}\theta(u)\ln u\,du.$$
Hence $g_\epsilon$ is bounded. By noticing that $\frac{r}{\epsilon^n}\leq \eta\Leftrightarrow n\leq \frac{\ln \frac{r}{\eta}}{\ln \epsilon}$, we deduce
\begin{align*}
\frac{\ln \frac{r}{\eta}}{\ln \epsilon}(k_\epsilon(0)-\alpha)+g_\epsilon(r)\leq K^\epsilon(r)\leq \frac{\ln \frac{r}{\eta}}{\ln \epsilon}k_\epsilon(0) +g_\epsilon(r).
\end{align*}
By taking the $\limsup$ and $\liminf$ in the above inequality, we have proved that for each $\alpha>0$:
$$ \frac{k_\epsilon(0)-\alpha}{\ln\frac{1}{\epsilon}} \leq  \liminf_{r\to 0}\frac{K^\epsilon(r)}{\ln\frac{1}{r}}\leq \limsup_{r\to 0}\frac{K^\epsilon(r)}{\ln\frac{1}{r}}\leq \frac{k_\epsilon(0)}{\ln\frac{1}{\epsilon}} ,$$
 which completes the proof.
\qed

\noindent {\it Proof of Corollary \ref{atom}.} By stationarity, it is enough to prove that, almost surely, the measure $M$ does not possess any atom on the segment $[0,1]$. 
From \cite[Corollary 9.3 VI]{daley}, it is enough to check that for each $\alpha>0$:
$$\sum_{k=1}^n\p\Big(M[\frac{k-1}{n};\frac{k}{n}]>\alpha\Big) = n\p\Big(M[\frac{0}{n};\frac{1}{n}]>\alpha\Big)\to 0\quad \text{ as }n\to \infty.$$
This is a direct consequence of the Markov inequality and Lemma \ref{estimate}:
$$n\p\Big(M[\frac{0}{n};\frac{1}{n}]>\alpha\Big)\leq \frac{n}{\alpha^{1+\gamma} }\E[M([0,\frac{1}{n}])^{1+\gamma}]\to 0\text{ as }n\to \infty.\qed$$

\noindent  {\it Proof of Proposition \ref{exactsi}.}
 Otherwise, if $M$ is a good lognormal $\star$-scale invariant random measure, then using Theorem \ref{main}, we know that there exists $k$ a 
 continuous covariance function such that, for all $|r|\leq R$: 
 \begin{equation}\label{eq_K}
 K(r) = \int_{|r|}^\infty \frac{k(u)}{u} du = \lambda^2\ln \frac{T}{|r|}+C. 
 \end{equation}
 By differentiating this equality with respect to $r$, we obtain $k(r)=\lambda^2$ for all $|r|\leq R$. Then, let $(X_t)_{t\in\R}$ be a centered stationary Gaussian process with covariance kernel $k$. For all $s,t \in \R$ such that $|t-s|<R$, we have $cov(X_t, X_s) = k(|t-s|) = k(0) = var[X_t]$
which implies (by Cauchy-Schwarz inequality) that $X_t =X_s $ almost surely. The process $X$ being stationary, this shows that it is a constant process.  
Hence $k(r)=\lambda^2$ for all $r \in \R$. Because of equation \eqref{eq_K}, this is a contradiction since it would imply $K(r)=+\infty$ for all $r$.
 \qed



\begin{thebibliography}{20}

\bibitem{applebaum}
{D. Applebaum}, \textit{L\'evy Processes and Stochastic Calculus}, Cambridge studies in advanced mathematics 93, Cambridge University Press, Cambridge, 2004.

\bibitem{bacry}
Bacry E., Muzy J.F.: Log-infinitely divisible multifractal processes, \emph{Comm. Math. Phys.}, \textbf{236} (2003) no.3, 449-475.


\bibitem{cf:BaKoMu} Bacry E., Kozhemyak, A., Muzy J.-F.: Continuous
  cascade models for asset returns, \emph{Journal of Economic Dynamics and Control}, \textbf{32} (2008) no.1, 156-199.


\bibitem{Bar}
Barral, J., Mandelbrot, B.B.: Multifractal products of cylindrical pulses, \emph{Probab. Theory
Relat. Fields} \textbf{124} (2002), 409-430.

\bibitem{Benj}
Benjamini, I., Schramm, O.: KPZ in one dimensional random geometry of multiplicative cascades, 

\bibitem{cf:Castaing}
Castaing B., Gagne Y., Hopfinger E.J.:
Velocity probability density-functions of high Reynolds-number turbulence,
\emph{Physica D} \textbf{46} (1990) 2, 177-200.

\bibitem{cf:Cas}
Castaing B., Gagne Y., Marchand M.:
Conditional velocity pdf in 3-D turbulence,
\emph{J. Phys. II France} \textbf{4} (1994), 1-8.



\bibitem{cf:DuRoVa} Duchon, J., Robert, R., Vargas, V.: Forecasting volatility with the multifractal random walk model, submitted to \emph{Mathematical Finance}, available at http://arxiv.org/abs/0801.4220.
\bibitem{daley} Daley D.J., Vere-Jones D., An introduction to the theory of point processes volume 2,  Probability and its applications, Springer, 2nd edition, 2007.


\bibitem{cf:DuSh} Duplantier, B., Sheffield, S.: Liouville Quantum Gravity and KPZ, to appear in \emph{Inventiones Mathematicae}, available on arxiv at the URL http://arxiv.org/abs/0808.1560.


\bibitem{cf:Fr} Frisch, U.: \emph{Turbulence}, Cambridge University Press (1995).

\bibitem{cf:Gel} Gelfand I.M., Vilenkin, N YA.: \emph{Generalized Functions}, Vol. 4, Academic Press, New 
York (1964).


\bibitem{cf:Kah} Kahane, J.-P.: Sur le chaos multiplicatif,
  \emph{Ann. Sci. Math. Qu{\'e}bec}, \textbf{9} no.2 (1985), 105-150.


\bibitem{cf: KPZ} Knizhnik, V.G., Polyakov, A.M., Zamolodchikov, A.B.: Fractal structure of 2D-quantum gravity, \emph{Modern Phys. Lett A}, \textbf{3}(8) (1988), 819-826.

\bibitem{mandelbrotstar}
Mandelbrot, B.B.: Multiplications al\'eatoires it\'er\'ees et distributions invariantes par moyenne pond\'er\'ee al\'eatoire, I and II.  \emph{Comptes Rendus} (Paris):  278A, 289-292 and 355-358. 


\bibitem{mandelbrot}
Mandelbrot B.B.: Intermittent turbulence in self-similar cascades, divergence of high moments and dimension of the carrier, \emph{J. Fluid. Mech.} \textbf{62} (1974), 331-358.




\bibitem{cf:RoVa} Robert, R., Vargas, V.: Gaussian Multiplicative Chaos revisited, \emph{Annals of Probability}, \textbf{38} 2 (2010), 605-631.

\bibitem{cf:RhoVar} Rhodes, R. Vargas, V.: KPZ formula for log-infinitely divisible multifractal random measures,  to appear in ESAIM, available on arxiv at the URL http://arxiv.org/abs/0807.1036. 

\bibitem{rhovar} Rhodes, R. Vargas, V.: Multidimensional multifractal random measures, \emph{Electronic Journal of Probability}, \textbf{15} (2010), 241-258.  

\bibitem{cf:Sch} Schmitt, F., Lavallee, D., Schertzer, D., Lovejoy, S.: Empirical determination of universal
multifractal exponents in turbulent velocity fields, \emph{Phys. Rev. Lett.} \textbf{68} (1992), 305-308.

\bibitem{cf:Schw} Schwartz, L.: \emph{Th\'eorie des distributions}, Hermann (1997).  



\bibitem{cf:Sto} Stolovitzky, G., Kailasnath, P., Sreenivasan, K.R.: Kolmogorov's Refined Similarity Hypotheses, \emph{Phys. Rev. Lett.} \textbf{69}(8) (1992), 1178-1181.  \end{thebibliography}
\end{document}